\documentclass[aos,MSNbibl,seceqn,citesort,dvips]{arximspdf}
\usepackage{mathbh}

\doi{10.1214/11-AOS965} 
\volume{40}
\issue{2}
\pubyear{2012}
\firstpage{832}
\lastpage{860}

\makeatletter

\newcommand{\cal}{\mathcal}

\newcommand{\R}{\mathbb{R}}
\newcommand{\N}{\mathbb{N}}
\newcommand{\E}{\mathbb{E}}
\newcommand{\Pro}{\mathbb P}

\newcommand{\one}{{\mathbh{1}}}

\newcommand{\cD}{{\cal D}}
\newcommand{\cF}{{\cal F}}
\newcommand{\cL}{{\cal L}}
\newcommand{\cS}{{\cal S}}
\newcommand{\cX}{{\cal X}}
\newcommand{\cZ}{{\cal Z}}

\newcommand{\RERM}{\hat f^{\mathrm{RERM}}_n}
\newcommand{\ERM}{\hat f^{\mathrm{ERM}}_n}

\newcommand{\argmin}{\mathop{\arg\min}}

\newtheorem{Theo}{Theorem}

\newtheorem{Theorem}{Theorem}[section]
\newtheorem{Lemma}[Theorem]{Lemma}

\newproclaim{Definition}[Theorem]{Definition}

\newtheorem{Proposition}[Theorem]{Proposition}

\newproclaim{Remark}[Theorem]{Remark}
\newproclaim{Assumption}{Assumption}[section]

\makeatother

\begin{document}
\begin{frontmatter}

\title{General nonexact oracle inequalities for classes with a
subexponential envelope}
\runtitle{Nonexact oracle inequalities}

\begin{aug}
\author[A]{\fnms{Guillaume} \snm{Lecu{\'e}}\corref{}\thanksref{t1}\ead[label=e1]{Guillaume.Lecue@univ-mlv.fr}}
\and
\author[B]{\fnms{Shahar} \snm{Mendelson}\thanksref{t2}\ead[label=e2]{shahar@tx.technion.ac.il}}
\runauthor{G. Lecu{\'e} and S. Mendelson}
\affiliation{CNRS, Universit{\'e} Paris-Est Marne-la-vall{\'e}e and
Technion, Israel~Institute~of~Technology}
\address[A]{CNRS\\
LAMA\\
Universit{\'e} Paris-Est Marne-la-vall{\'e}e, 77454\\
France\\
\printead{e1}} 
\address[B]{Department of Mathematics\\
Technion, Israel Institute of Technology\\
Haifa 32000\\
Israel\\
\printead{e2}}
\end{aug}

\thankstext{t1}{Supported by the French Agence Nationale de
la Recherce (ANR) Grant ``\textsc{Prognostic}'' ANR-09-JCJC-0101-01.}

\thankstext{t2}{Supported in part by the Mathematical Sciences Institute---The
Australian National University, the European Research Council (under
ERC Grant agreement 203134) and the Australian Research Council
(under Grant DP0986563).}

\received{\smonth{1} \syear{2011}}
\revised{\smonth{12} \syear{2011}}

%
\begin{abstract}
We show that empirical risk minimization procedures and regularized
empirical risk minimization procedures satisfy nonexact oracle
inequalities in an unbounded framework, under the assumption that the
class has a subexponential envelope function. The main novelty, in
addition to the boundedness assumption free setup, is that those
inequalities can yield fast rates even in situations in which exact
oracle inequalities only hold with slower rates.

We apply these results to show that procedures based on $\ell_1$ and
nuclear norms regularization functions satisfy oracle inequalities with
a residual term that decreases like $1/n$ for every $L_q$-loss
functions (\mbox{$q\geq2$}), while only assuming that the tail behavior of the
input and output variables are well behaved. In particular, no RIP type
of assumption or ``incoherence condition'' are needed to obtain fast
residual terms in those setups. We also apply these results to the
problems of convex aggregation and model selection.
\end{abstract}

%
\begin{keyword}[class=AMS]
\kwd[Primary ]{62G05}
\kwd[; secondary ]{62H30}
\kwd{68T10}.
\end{keyword}
\begin{keyword}
\kwd{Statistical learning}
\kwd{fast rates of convergence}
\kwd{oracle inequalities}
\kwd{regularization}
\kwd{classification}
\kwd{aggregation}
\kwd{model selection}
\kwd{high-dimensional data}.
\end{keyword}

\end{frontmatter}

\section{Introduction and main results}
\label{secintro}
Let $\cZ$ be a space endowed with a probability measure $P$, and let
$Z$ and $Z_1,\ldots,Z_n$ be $n+1$ independent random variables with
values in
$\cZ$, distributed according to $P$; from the statistical point of view,
$\cD=(Z_1,\ldots,Z_n)$ is the set of given data. Let $\ell$ be a loss
function which associates a real number $\ell(f,z)$ to any real-valued
measurable function $f$ defined on $\cZ$ and any point $z\in\cZ$.
Denote by $\ell_f$ the loss function $\ell(f,\cdot)$ associated with
$f$ and set $R(f)=\E\ell_f(Z)$ to be the associated risk. The risk of
any statistic $\hat f_n(\cdot)=\hat f_n(\cdot,\cD)\dvtx\cZ
\longrightarrow\R
$ is defined by $R(\hat f_n)=\E[\ell_{\hat f_n}(Z)|\cD]$.

Let $F$ be a class (usually called the \textit{model}) of real-valued
measurable functions defined on $\cZ$. In learning theory, one wants to
assume as little as possible on the class $F$, or on the measure $P$.
The aim is to use the data to construct learning algorithms whose risk
is as close as possible to $\inf_{f\in F}R(f)$ (and when this infimum
is attained by a function $f^*_F$ in $F$, this element is called an
\textit{oracle}). Hence, one would like to construct procedures~$\hat
f_n$ such that, for some $\epsilon\geq0$, with high probability,
%
\begin{equation}
\label{eqoracle-inequality}
R(\hat f_n)\leq(1+\epsilon)\inf_{f\in F}R(f)+r_n(F).
\end{equation}
The role of the \textit{residual term} (or \textit{rate}) $r_n(F)$ is
to capture the ``complexity'' of the problem, and the hope is to make
it as small as possible.

When $r_n(F)$ tends to zero as $n$ tends to infinity, inequality (\ref
{eqoracle-inequality}) is called an \textit{oracle inequality}. When
$\epsilon=0$, we say that $\hat f_n$ satisfies an \textit{exact oracle
inequality} (the term \textit{sharp} oracle inequality has been also
used) and when $\epsilon>0$ it satisfies a \textit{nonexact oracle
inequality}. Note that the terminology ``risk bounds'' has been also
used for (\ref{eqoracle-inequality}) in the literature.

A natural algorithm in this setup is the \textit{empirical risk
minimization procedure} (ERM) (terminology due to \cite{MR672244}), in
which the \textit{empirical risk} functional
\[
f\longmapsto R_n(f)=\frac{1}{n}\sum_{i=1}^n\ell_f(Z_i)
\]
is minimized and produces $\hat f_n^{\mathrm{ERM}}\in\operatorname{Arg}\min_{f\in
F}R_n(f)$. Note that when $R_n(\cdot)$ does not achieve its infimum
over $F$ or if the minimizer is not unique, we define $\hat f_n^{\mathrm{ERM}}$
to be an element in $F$ for which $R(\hat f_n^{\mathrm{ERM}})\leq\inf_{f\in
F}R(f)+1/n$. This algorithm has been extensively studied, and we will
compare our first result to the one of \mbox{\cite{MR2329442,MR2240689,MR2291502}}.

One motivation in obtaining nonexact oracle inequalities
[equation (\ref{eqoracle-inequality}) for $\epsilon\not=0$] is the
observation that in many situations, one can obtain such an inequality
for the ERM procedure with a residual term $r_n(F)$ of the order of
$1/n$, while the best residual term achievable by ERM in an exact
oracle inequality [equation (\ref{eqoracle-inequality}) for $\epsilon
=0$] will only be of the order of $1/\sqrt{n}$ for the same problem.
For example, consider the simple case of a finite model $F$ of
cardinality $M$ and the bounded regression model with the quadratic
loss function [i.e., $Z=(X,Y)\in\cX\times\R$ with $|Y|,{\max_{f\in
F}}|f(X)|\leq C$ for some absolute constant $C$ and $\ell
(f,(X,Y))=(Y-f(X))^2$]. It can be verified that for every $x>0$, with
probability greater than $1-8\exp(-x)$, $\hat f_n^{\mathrm{ERM}}$ satisfies a
nonexact oracle inequality with a residual term proportional to
$(x+\log M)/(\epsilon n)$. On the other hand, it is known \cite
{lbw96,MR2451042,LM2} that in the same setup, there are finite models
for which, with probability greater than a positive constant, $\hat
f_n^{\mathrm{ERM}}$ cannot satisfy an exact oracle inequality with a residual
term better than $c_0\sqrt{(\log M)/n}$. Thus, it is possible to
establish two optimal oracle inequalities [i.e., oracle inequalities
with a nonimprovable residual term $r_n(F)$ up to some multiplying
constant] for the same procedure with two very different residual
terms: one being the square of the other one. We will see below that
the same phenomenon occurs in the classification framework for VC
classes. Thus our main goal here is to present a general framework for
nonexact oracle inequalities for ERM and RERM (regularized ERM), and
show that they lead to fast rates in cases when the best known exact
oracle inequalities have slow rates.

Although the improved rates are significant, it is clear that exact
inequalities are more ``valuable'' from the statistical point of view.
For example, consider the regression model with the quadratic loss. It
follows from an exact oracle inequality on the prediction risk
[equation (\ref{eqoracle-inequality}) for $\epsilon=0$], another exact
oracle inequality, but for the estimation risk
\[
\|\hat f^{\mathrm{ERM}}_n-f^*\|_{L_2}^2\leq\inf_{f\in F}\|f-f^*\|_{L_2}^2+r_n(F),
\]
where $f^*$ is the regression function of $Y$ given $X$, and $\|\cdot\|
_{L_2}$ is the $L_2$-norm with respect to the marginal
distribution of $X$.

In other words, exact oracle inequalities for the prediction risk
$R(\cdot)$ provide both prediction and estimation results (prediction
of the output $Y$ and estimation of the regression function $f^*$)
whereas nonexact oracle inequalities provide only prediction results.

Of course, nonexact inequalities are very useful when it suffices
to compare the risk $R(\hat f_n)$ with $(1+\epsilon)\inf_{f\in F}R(f)$;
and the aim of this note is to show that the residual term can be
dramatically improved in such cases.
\subsection{Empirical risk minimization}
\label{secERM}
The first result of this note is a nonexact oracle inequality for the
ERM procedure. To state this result, we need the following notation.
Let $G$ be a class of real-valued functions defined on~$\cZ$. An
important part of our analysis relies on the behavior of the supremum
of the empirical process indexed by $G$
%
\begin{equation}\label{eqsupremum-empirical-process}
\|P-P_n\|_G={\sup_{g\in G}}|(P-P_n)(g)|,
\end{equation}
where for every $g\in G$, we set $P g= \E g(Z)$ and $P_n g=n^{-1}\sum
_{i=1}^ng(Z_i)$. Recall that for every $\alpha\geq1$, the $\psi
_\alpha$
norm of $g(Z)$ is
\[
\|g(Z)\|_{\psi_\alpha}=\inf\bigl(c>0\dvtx\E\exp\bigl(|g(Z)|^\alpha
/c^\alpha
\bigr)\leq2\bigr).
\]
We will control the supremum (\ref{eqsupremum-empirical-process})
using the quantities
\[
\sigma(G)=\sup_{g\in G}\sqrt{P g^2} \quad\mbox{and}\quad b_n(G)=\Bigl\|{\max
_{1\leq i\leq n}\sup_{g\in G}}|g(Z_i)|\Bigr\|_{\psi_1}.
\]
Note that for a bounded class $G$, one has $b_n(G)\leq\sup_{g\in
G}\|g\|_{\infty}$ and in the sub-exponential case, $b_n(G)\lesssim
(\log en)\|{\sup_{g\in G}}|g|\|_{\psi_1}$ (this follows from Pisier's
inequality); cf. Lemma 2.2.2 in \cite{vanderVaartWellner}.
Throughout this note we will also use the notation $b_n(g)=\|{\max
_{1\leq i\leq n}}|g(Z_i)|\|_{\psi_1}$ and for any pseudo-norm~$\|\cdot
\|$ on $L_2(P)$, we will denote by $\operatorname{diam}(G,\|\cdot\|)={\sup
_{g\in G}}\|g\|$ the diameter of~$G$ with respect to this norm.

Observe that the desired bound depends on the $\psi_1$ behavior of the
\textit{envelope function} of the class, ${\sup_{g \in G}}
|g(Z)|$,\vspace*{1pt} and as noted above, this extends the
``classical'' framework of a uniformly bounded class in $L_\infty$.
Although this extension seems minor at first, the examples we will
present show that the assumption is not very restrictive and allows one
to deal with LASSO-type situations, in which the indexing class is very
small---something which is impossible under the $L_\infty$ assumption.
On the other hand, it should be emphasized that this is not a step
towards an unbounded learning theory. For such results, the analogous
assumption should be that the class has a bounded diameter in $\psi_1$,
which is, of course, a much weaker assumption than a $\psi_1$ envelope
function and requires different methods; see, for example, \cite{Shahar-Gaussian,Shahar-Paouris}.

To obtain the required bound, we will study empirical processes indexed
by sets associated with $G$, namely, the star-shaped hull of $G$ around
zero and the localized subsets for different levels $\lambda\geq0$,
defined by
\[
V(G)=\{\theta g\dvtx 0\leq\theta\leq1, g\in G\} \quad\mbox{and}\quad
V(G)_\lambda
=\{h\in V(G)\dvtx Ph\leq\lambda\}.
\]
%

Given a model $F$ and a loss function $\ell$, consider the loss class
and the excess loss class $ \ell_F\,{=}\,\{\ell_f\dvtx f\,{\in}\,F\}$ and the excess
loss class \mbox{$\cL_F\,{=}\,\{\ell_f\,{-}\,\ell_{f^*_F}\dvtx f\,{\in}\, F\}$}. We will assume that
an oracle $f^*_F$ exists in $F$, and from here on set $\cL_f=\ell
_f-\ell
_{f^*_F}$.
\renewcommand{\theTheo}{\Alph{Theo}}
\begin{Theo}\label{TheoA}
There exists an absolute constant
$c_0>0$ for which the following holds. Let $F$ be a class of functions
and assume that there exists $ B_n\geq0$ such that for every $f\in F$,
$P\ell_f^2\leq B_n P\ell_f+B_n^2/n$. Let $0<\epsilon<1/2$, set
$\lambda
^*_\epsilon>0$ for which
\[
\E\|P_n - P\|_{V(\ell_F)_{\lambda^*_\epsilon}} \leq(\epsilon
/4)\lambda
^*_\epsilon
\]
and put $\rho_n$ an increasing function satisfying that for every $x>0$,
\[
\rho_n(x)\geq\max\biggl(\lambda^*_\epsilon,c_0\frac{(b_n(\ell
_F)+B_n/\epsilon)x}{n\epsilon}\biggr).
\]
Then, for every $x>0$, with probability greater than $1-8\exp(-x)$,
\[
R(\hat f_n^{\mathrm{ERM}})\leq(1+3\epsilon)\inf_{f\in F}R(f)+\rho_n(x).
\]
\end{Theo}
%
\begin{Remark}
Although the formulation of Theorem \ref{TheoA} requires that for every $\ell\in
\ell_F$, $P\ell^2\leq B_n P\ell+B_n^2/n$, we will show that if $\ell$
is nonnegative, this condition is trivially satisfied for $B_n \sim
\operatorname{diam}(\ell_F,\psi_1)\log(n)$.\vadjust{\goodbreak}

Unfortunately, this type of condition is far from being trivially
satisfied for the excess loss class $\cL_F=\{\ell_f-\ell
_{f^*_F}\dvtx f\in
F\}$, which is one of the major differences between exact and nonexact
oracle inequalities. Indeed, the \textit{Bernstein condition}, that for
every $f \in F$, $\E\cL_f^2\leq B\E\cL_f$ (see\vspace*{1pt} \cite{MR2240689} or
Section \ref{secComments} below), used in \mbox{\cite
{MR2291502,MR2329442,MR2240689}} to obtain exact oracle inequalities
with fast rates (rates of the order of $1/n$), depends on the geometry
of the problem \cite{MR2426759,m08} and may not be true in general.
Theorem \ref{TheoA} is similar in nature to Corollary~2.9 of \cite{MR2240689} and
a detailed comparison between the two results can be found in
Section \ref{secComments}.
\end{Remark}

Theorem \ref{TheoA} is similar in nature to Theorem 2 in
\cite{MR2291502}.
%
\begin{Theorem}
Let $\phi\dvtx \R\to\R$ be a nondecreasing, continuous function, for
which $\phi(1)\geq1$ and $x\rightarrow\phi(x)/x$ is nonincreasing.
Set $F$ to be a class of functions where there is some $0\leq\beta
\leq
1$ such that $\E\cL_f^2\leq B(\E\cL_f)^{\beta}$ and
$\|\ell _f\|_\infty\leq1$.
If $\phi(\lambda)\geq\sqrt{n}\E\sup_{f,g\in F, P(\ell_f-\ell
_g)^2\leq
\lambda^2 }(P-P_n)(\ell_f-\ell_g)$ for any $\lambda$ satisfying
$\phi
(\lambda)\leq\sqrt{n}\lambda^2$, and $\varepsilon_*$ is the unique
solution of the equation $\sqrt{n} \varepsilon_*^2=\phi(\sqrt
{B}\varepsilon_*^{\beta})$, then for every $x\geq1$, with probability
greater than $1-\exp(-x)$,
\[
R(\hat f_n^{\mathrm{ERM}})\leq\inf_{f\in F}R(f)+c_0 x \varepsilon_*^2.
\]
\end{Theorem}

One of the applications of the above theorem in learning theory is for
the loss function $\ell_f(x,y)=\one_{f(x)\neq y}$. It leads to an exact
oracle inequality for the ERM procedure, preformed in a class $F$ of VC
dimension \mbox{$V\leq n$} (see \cite{MR2291502} for more details), and with a
residual term of the order of $(V\log(enB^{1/\beta}/V)/\break n
)^{1/(2-\beta)}$.

In comparison, in the same situation, for every $f\in F$, $\E\ell
_f^2\leq\E\ell_f$. Therefore, it follows from Theorem \ref{TheoA}, the argument
used to obtain equation (29) in \cite{MR2291502} (or Example 3 in
\cite
{MR2329442}) and the peeling argument which will be presented in (\ref
{eqpeeling-shahar}) below, that for every $x\geq1$, with probability
greater than $1-8\exp(-x)$,
%
\begin{equation}\label{eqMassart-NEOI}
R(\hat f_n^{\mathrm{ERM}})\leq(1+3\epsilon)\inf_{f\in F}R(f)+c_0\frac{xV\log
(en/V)}{\epsilon^2n}.
\end{equation}

The residual term $\epsilon_*^2$ obtained in \cite{MR2291502} is
optimal, but since it heavily depends on the parameter $\beta$, it
ranges between $\sqrt{V/n}$ and $V/n$ (up to a~logarithmic factor). In
particular, it can be as bad as the \textit{square root} of the
residual term of the nonexact oracle inequality (\ref
{eqMassart-NEOI}) in the same situation. The main difference between
the two results is that the condition $\E\ell_f^2\leq\E\ell_f$ for
every $f \in F$ is always satisfied whereas the condition that for
every $f\in F$ $ \E\cL_f^2\leq B (\E\cL_f)^{\beta}$
depends on
the relative position of $Y$ and $F$, and thus on geometry of the
system $(F,Y)$.

It is interesting to note that the residual term in (\ref
{eqMassart-NEOI}) always yields fast rate even for hard classification\vadjust{\goodbreak}
problem such that $\Pro[Y=1|X]=1/2$. This means that while the
prediction problem in classification is completely blind to the
geometry of the model, the estimation problem is influenced in a~very
strong way by the geometry of $(F,Y)$. Thus, estimating the regression
function (or the Bayes rule) is in general much harder than predicting
the output~$Y$.

Another related result is the one in \cite{MR2329442} where (among
other results) an exact oracle inequality is proved for the ERM with a
residual term~$\delta_n(x)$. The residual term is controlled using the empirical
oscillation\vspace*{1pt} $\phi_n(\delta)\,{=}\,{\E\sup_{f,g\in F(\delta)}}|\allowbreak (P-P_n)(\ell
_f-\ell_g)|$ indexed by $F(\delta)=\{f\in F\dvtx P\cL_f\leq\delta\}$, and
by the $L_2$ diameter $D(\delta)=\sup_{f,g\in F(\delta)}\sqrt
{P(\ell
_f-\ell_g)^2}$
\[
\delta_n(x)=\argmin\Biggl(\delta>0\dvtx \phi_n(\delta) +\sqrt{\frac
{2x}{n}
\bigl(D(\delta)^2+2\phi_n(\delta)\bigr)}+\frac{x}{2n}\leq c_0 \delta
\Biggr).
\]
%

Note that all the quantities $\lambda_\epsilon^*$, $\varepsilon_*^2$
from \cite{MR2291502}, $\delta_n(x)$ from \cite{MR2329442}, $\mu^*$
from \cite{MR2240689} or Theorem \ref
{theoExact-oracle-ineq-ERM} below, define the residual terms of the
oracle inequalities as a fixed point of some equation. Those appear
naturally either from \textit{iterative localization} of the excess
risk, converging to $\delta_n(x)$ \cite{MR2329442,MR1857339}, or from
an ``isomorphic'' argument \cite{MR2240689} identifying the ``level''
$\mu
^*$ at which the actual and the empirical structures are equivalent. We
refer the reader to those articles for more details.

Results in \cite{MR2291502,MR2329442,MR2240689} were obtained under the
boundedness assumption ${\sup_{f\in F}}\|\ell_f\|_\infty\leq1$ because
the necessary tools from empirical processes theory, like contraction
inequalities \cite{LT91}, only hold under such an assumption. In
particular, these results do not apply even to the Gaussian regression
model. The approach developed in this work provides a slight
improvement, since risk bounds hold if the envelope function $\sup
_{f\in F}\ell_f$ is sub-exponential (which is the case for the Gaussian
regression model with respect to the square loss).


One should also mention the subtle but significant gap between the
\textit{margin assumption} and the Bernstein condition which we use. Both state
that for every $f \in F$,
\[
\E(\ell_f-\ell_{f^*})^2\leq B_0\bigl(\E(\ell_f-\ell_{f^*})
\bigr)^{1/\kappa}
\]
for some constant $\kappa\geq1$. However, in the margin condition
$f^*$ has the minimal risk \textit{over all measurable functions} (for
instance, $f^*$ is the regression function in the regression model with
respect to the quadratic loss), while in a Bernstein condition $f^*_F$
is assumed to minimize the risk \textit{over $F$}.

The two conditions are equivalent only when $f^*\in F$ (and thus
$f^*=f^*_F$). But in general, they are very different. As a simple
example, in the bounded regression model [i.e., $|Y|,{\sup_{f\in
F}}|f(X)|\leq C$] with respect to the quadratic loss, the margin
assumption holds with $\kappa=1$ whereas the Bernstein condition is not
true in general. For more details on the difference between the margin
assumption and the Bernstein condition we refer the reader to the
discussion in \cite{LM4}.

\subsection{Regularized empirical risk minimization}
\label{secRERM}
The second type of application we will present deals with nonexact
regularized oracle inequalities. Usually a model $F$ is chosen or
constructed according to the belief that an oracle $f^*_F$ in $F$ is
close, in some sense, to some minimizer $f^*$ of the risk function in
some larger class of functions~$\cF$ [e.g., in the regression model,~$f^*$ can be the regression function and $\cF=L^2(P_X)$]. Hence, by
choosing a~particular model $F\subset\cF$, it implicitly means that we
believe $f^*$ to be close to~$F$ in some sense.

It is not always possible to construct a class $F$ that captures
properties~$f^*$ is believed to have (e.g., a low-dimensional structure
or some smoothness properties). In such situations, one is not given a
single model $F$ (usually the set $\cF$ is too large to be called a
model), but a functional $\operatorname{crit}\dvtx\cF\longrightarrow\R^+$, called a
\textit{criterion}, that characterizes each function according to its
level of compliance with the desired property---and the smaller the
criterion, the ``closer'' one is to the property. For instance, when
$\cF$ is an RKHS, one can take $\operatorname{crit}(\cdot)$ to be the norm in
the reproducing kernel Hilbert space, or when~$\cF$ is the set of all
linear functionals in $\R^d$, one may chose $\operatorname{crit}(\beta)=\|
\beta\|_{\ell_p}$ for some $p\in[0,\infty]$. The extreme case here is
$p=0$ and $\|\beta\|_{\ell_0}$ is the cardinality of the support of
$\beta$; thus a small criterion means that $\beta$ belongs to a
low-dimensional space.

Instead of considering the ERM over the too large class $\cF$, the goal
is to construct a procedure having both good empirical performances and
a~small criterion. One idea, that we will not develop here, is to
minimize the empirical risk over the set
$F_r=\{f\in\cF\dvtx\operatorname{crit}(f)\leq r\}$
\cite{MR1379242,BMN}, and try to find a~data-dependent way of choosing
the radius $r$. Another popular idea is to regularize the empirical
risk: consider a nondecreasing function of the criterion called a
\textit{regularizing function} and denoted by
$\operatorname{reg}\dvtx\cF\longrightarrow\R^+$ and construct
%
\begin{equation}
\label{eqRERM} \hat
f_n^{\mathrm{RERM}}\in\mathop{\operatorname{Arg}\min}_{f\in\cF}\bigl(R_n(f)+\operatorname{reg}(f)\bigr)
\end{equation}
with the obvious extension if the infimum is not attained.

The procedure (\ref{eqRERM}) is called \textit{regularized empirical
risk minimization procedure} (RERM). RERM procedures were introduced to
avoid the ``over-fitting'' effect of large models
\cite{MR1679028,MR2319879}, and later used to select functions with
additional properties, like smoothness (e.g., SVM estimators in
\cite{MR2450103}) or an underlying low-dimensional structure (e.g., the
LASSO estimator).

In this setup, we are interested in constructing estimators $\hat f_n$
realizing the best possible trade-off between the risk and the
regularizing function over~$\cF$: there exists some $\epsilon\geq0$
such that with high probability
%
\begin{equation}
\label{eqOracle-Ineq-penalized}
R(\hat f_n)+\operatorname{reg}(\hat f_n)\leq(1+\epsilon)\inf_{f\in\cF}
\bigl(R(f)+\operatorname{reg}(f)\bigr).\vadjust{\goodbreak}
\end{equation}
Using the same terminology as in (\ref{eqoracle-inequality}),
inequality (\ref{eqOracle-Ineq-penalized}) is called a
\textit{regularized oracle inequality}. When $\epsilon=0$, (\ref
{eqOracle-Ineq-penalized}) is called an \textit{exact regularized
oracle inequality}, and when $\epsilon>0$, (\ref
{eqOracle-Ineq-penalized}) is called a \textit{nonexact regularized
oracle inequality}.

Following our analysis of the ERM algorithm, the next result is a
regularized oracle inequality for the RERM. But before stating this
result, one has to say a word on the way the regularizing function
$\operatorname{reg}(\cdot)$ and the criterion $\operatorname{crit}(\cdot)$ are related.

The choice of $\operatorname{reg}(\cdot)$ is driven by the complexity of the
sequence $(F_r)_{r\geq0}$ of models
\[
F_r=\{f\in\cF\dvtx\operatorname{crit}(f)\leq r\}.
\]
For any $r\geq0$, the complexity of $F_r$ is measured by $\lambda
_\epsilon^*(r)$ defined as above for some fixed $0<\epsilon<1/2$ by
\[
\E\|P_n - P\|_{V(\ell_{F_r})_{\lambda_\epsilon^*(r)}} \leq
(\epsilon
/4)\lambda_\epsilon^*(r).
\]
Hence, $\lambda_\epsilon^*(r)$ is a ``level'' in $\ell_{F_r}$ above
which the empirical and the actual structures are equivalent; namely,
with high probability, on the set $\{\ell\in\ell_{F_r}\dvtx\allowbreak P\ell
\geq
\lambda_\epsilon^*(r)\}$,
\[
(1/2)P_n\ell\leq P\ell\leq(3/2)P_n\ell.
\]
Thus, the function $r\rightarrow\lambda_\epsilon^*(r)$ captures the
``isomorphic profile'' of the collection $(\ell_{F_r})_{r \geq0}$. Up
to minor technical adjustments, the regularizing function, defined
formally in (\ref{eqpenalty-tilde}), is
$\operatorname{reg}(\cdot)=\lambda_\epsilon^*(\operatorname{crit}(\cdot))$.

We will study two separate situations, both motivated by the
applications we have in mind. In the first, $\operatorname{crit}(\cdot)$ will be
uniformly bounded and may only grow with the sample size $n$---that
is, there is a constant $C_n$ satisfying that for every \mbox{$f \in\cF$},
$\operatorname{crit}(f)\leq C_n$. The second case we deal with is when the
``isomorphic profile'' $r\rightarrow\lambda_\epsilon^*(r)$ tends to
infinity with $r$. For technical reasons, we also introduce an
auxiliary function $\alpha_n$, defined in the following assumption.
%
\begin{Assumption} \label{assthm-B}
Assume that for every $f\in\cF$, $\ell_f(Z)\geq0$ a.s. and that there
are nondecreasing functions $\phi_n$ and $B_n$ such that for every
$r\geq0$ and every $f \in F_r$,
\[
b_n(\ell_{F_r})\leq\phi_n(r) \quad\mbox{and}\quad P\ell_f^2\leq B_n(r)
P\ell
_f+B_n^2(r)/n.
\]
Let $0<\epsilon<1/2$ and consider a function $\rho_n\dvtx\R_+\times\R_+^*
\to\R$ nondecreasing in its first argument and such that, for any
$r\geq0$ and $x>0$,
\[
\rho_n(r,x)\geq\max\biggl(\lambda_\epsilon^*(r),c_0\frac{(\phi
_n(r)+B_n(r)/\epsilon)(x+1)}{n\epsilon}\biggr).
\]

Assume that either:
\begin{itemize}
\item there exists $C_n>0$ such that for every $f\in\cF,\operatorname{crit}(f)\leq C_n$ and in this case define $\alpha_n(\epsilon,x)= C_n$,
for all $0<\epsilon<1/2$ and $x>0$, or\vadjust{\goodbreak}
\item the function $r\rightarrow\lambda_\epsilon^*(r)$ tends to
infinity with $r$ and there exists $K_1>0$ such that $2\rho_n(r,x)\leq
\rho_n(K_1(r+1),x)$, for all $r\geq0$ and $x>0$ and, in this case, let
$f_0$ be any function in $\bigcup_{r\geq0}F_r$ and define $\alpha_n$ such
that, for every $x>0$ and $0<\epsilon<1/2$,
%
\begin{eqnarray}
\label{eqalpha-function}\quad
&&\alpha_n(\epsilon,x)\geq\max\bigl[K_1\bigl(\operatorname{crit}(f_0)+2\bigr),\nonumber\\
&&\hspace*{72.2pt}(\lambda_\epsilon^*)^{-1}\bigl((1+2\epsilon)\bigl(3R(f_0)+2K^\prime
\bigl(b_n(\ell
_{f_0})+B_n(\operatorname{crit}(f_0))\bigr)\\
&&\hspace*{253pt}{}\times\bigl((x+1)/n\bigr)\bigr)\bigr)\bigr],\nonumber
\end{eqnarray}
where $(\lambda_\epsilon^*)^{-1}$ is the generalized inverse function
of $\lambda_\epsilon^*$ [i.e., $(\lambda_\epsilon^*)^{-1}(y)=\sup
(r>0\dvtx \lambda_\epsilon^*(r)\leq y)$, for all $y>0$] and
$K^\prime$
is some absolute constant.
\end{itemize}
\end{Assumption}
\begin{Theo}\label{TheoB}
There exist absolute positive constants $c_0$, $c_1$ $K$ and $K^\prime$
for which the following holds. Under Assumption \ref{assthm-B}, for
every $x>0$ and
%
\begin{equation}\label{eqRERM-theo-peter}
\hat f_n^{\mathrm{RERM}}\in\mathop{\operatorname{Arg}\min}_{f\in\cF}
\biggl(R_n(f)+\frac
{2}{1+2\epsilon}\rho_n\bigl(\operatorname{crit}(f)+1,x+\log\alpha_n(\epsilon
,x)\bigr)\biggr),\hspace*{-35pt}
\end{equation}
with probability greater than $1-12\exp(-x)$,
\begin{eqnarray*}
&& R(\hat f_n^{\mathrm{RERM}})+\rho_n\bigl(\operatorname{crit}(\hat f_n^{\mathrm{RERM}})+1,x+\log
\alpha
_n(\epsilon,x)\bigr)\\
&&\qquad\leq\inf_{f\in\cF}\biggl[(1+3\epsilon)R(f)+2\rho_n\bigl(\operatorname{crit}(f)+1,x+\log\alpha_n(\epsilon,x)\bigr)\\
&&\qquad\quad\hspace*{76pt}{} +c_1\frac{(b_n(\ell_f)+B_n(\operatorname{crit}(f))/\epsilon
)(x+1)}{n\epsilon}\biggr].
\end{eqnarray*}
\end{Theo}

Fortunately, $\alpha_n$ usually has little impact on the resulting
rates. For instance, in the main application we will present here,
\mbox{$\log\alpha_n(\epsilon,x)\lesssim_\epsilon\log(x+n)$}.

Like in Theorem \ref{TheoA}, the Bernstein-type condition $ P\ell^2\leq B_n(r)
P\ell+B_n^2(r)/n$ holds when $\ell$ is nonnegative and sub-exponential
for $B_n(r)\lesssim\operatorname{diam}(\ell_{F_r},\psi_1)\times\allowbreak \log(n)$. Therefore,
and contrary to the situation in exact oracle inequalities, the
``geometry'' of the family of classes $(F_r)_{r\geq0}$ does not play a
crucial role in the resulting nonexact regularized oracle inequalities.

Observe that now the choice of the regularizing function in terms of
the criterion is now made explicit:
%
\begin{equation}\label{eqpenalty-tilde}
\operatorname{reg}(f)=\frac{2}{1+2\epsilon}\rho_n\bigl(\operatorname{crit}(f)+1,x+\log
\alpha
_n(\epsilon,x)\bigr).
\end{equation}



\subsection{\texorpdfstring{$\ell_1$-regularization}{l_1-regularization}}
\label{secl1-regularization}
The formulation of Theorem \ref{TheoB} seems cumbersome, but it is not very
difficult to apply it---and here we will present one application
dealing with high-dimensional vectors of short support. Other
applications on matrix completion, convex aggregation and model
selection can be found in \cite{supplementary}.

Formally, let $(X,Y), (X_i,Y_i)_{1\leq i\leq n}$ be $n+1$ i.i.d. random
variables with values in $\R^d\times\R$, and denote by $P_X$ the\vadjust{\goodbreak}
marginal distribution of $X$. The dimension $d$ can be much larger than
$n$ but we believe that the output $Y$ can be well predicted by a
sparse linear combination of covariables of $X$; in other words, $Y$
can be reasonably approximated by $\langle X,\beta_0 \rangle$ for
some $\beta
_0\in
\R^d$ of short support (even though we will not require any assumption
of this type to obtain our results).

These kind of problems are called ``high-dimensional'' because there are
more covariables than observations. Nevertheless, one hopes that under
the structural assumption that $Y$ ``depends'' only on a few number of
covariables of $X$, it would still be possible to construct efficient
statistical procedures to predict $Y$.

In this framework, a natural criterion function is the $\ell_0$
function measuring the size of the support of a vector. But since this
function is far from being convex, using it in practice is hard; see,
for example, \cite{MR1320206}. Therefore, it is natural to consider a
convex relaxation of the~$\ell_0$ function as a criterion: the~$\ell_1$
norm \cite{MR1379242,MR2382644,MR2241189}.

In what follows, we will apply Theorem \ref{TheoB} to establish nonexact
regularized oracle inequalities for $\ell_1$-based RERM procedures, and
with fast error rates---a~residual term that tends to $0$ like $1/n$
up to logarithmic terms. The regularizing function resulting from
Theorem \ref{TheoB} for the $L_q$-loss ($q\geq2$) will be the $q$th power of the
$\ell_1$-norm. In particular, for the quadratic loss, we regularize by
$\|\cdot\|^2_{\ell_1}$, the \textit{square of the $\ell_1$-norm},
%
\begin{equation}\label{eqSquare-LASSO}
\hat\beta_n\in\mathop{\operatorname{Arg}\min}_{\beta\in\R^d}\Biggl(\frac{1}{n}\sum
_{i=1}^n(Y_i-\langle X_i,\beta\rangle)^2
+\kappa(n,d,x)\frac{\|\beta\|_{\ell_1}^2}{n}\Biggr),
\end{equation}
while the standard LASSO is regularized by the $\ell_1$ norm itself.
This choice of the exponent is dictated by the complexity of the
underlying models: the sequence of balls $(rB_1^d)_{r\geq0}$ trough the
isomorphic profile function $r\rightarrow\lambda_\epsilon^*(r)$.
Observe that since $\|\beta\|_{\ell_1}/\sqrt{n}\geq\|\beta \|
_{\ell_1}^2/n$ when $\|\beta\|_{\ell_1}\leq\sqrt{n}$, a nonexact
oracle inequality for the LASSO estimator itself follows from Theorem
\ref{TheoB}, but with a slow rate of $1/\sqrt{n}$. Using the $q$th power of the
$\ell_1$-norm as a penalty function for the $L_q$-risk yields a fast
$1/n$ rate (see Theorem \ref{TheoC}). 


We will perform this study for the $L_q$-loss function, and in which
case, for every $\beta\in\R^d$,
\[
R^{(q)}(\beta)=\E|Y-\langle X,\beta\rangle|^q
\quad\mbox{and}\quad
R^{(q)}_n(\beta
)=\frac
{1}{n}\sum_{i=1}^n|Y_i-\langle X_i,\beta\rangle|^q.
\]
The following result is obtained only under the assumption that $Y$ and
$\|X\|_{\ell_\infty^d}$ belong to $L_{\psi_q}$. Since there are no
``statistically reasonable'' $\psi_q$ variables for \mbox{$q>2$}, it sounds
more ``statistically relevant'' to assume that $|Y|$, $\|X\|_{\ell
_\infty^d}$ are almost surely bounded when one wants results for the
$L_q$-risk with $q>2$, or that the functions are in $L_{\psi_2}$ for
$q=2$ (e.g., linear models with sub-Gaussian noise and a sub-Gaussian
design satisfy this condition).

\begin{Theo}\label{TheoC}
Let $q\geq2$. There exist constants
$c_0$ and $c_1$ that depend only on $q$ for which the following holds.
Assume that there exists $K(d)>0$ such that $\|Y\|_{\psi_q}$,
$\|\|X\|_{\ell_\infty^d}\|_{\psi_q}\leq K(d)$. For $x>0$ and
$0<\epsilon
<1/2$, let
\[
\lambda(n,d,x)=c_0 K(d)^q(\log n)^{(4q-2)/q}(\log d)^2(x+\log n)
\]
and consider the RERM estimator
\[
\hat\beta_n\in\mathop{\operatorname{Arg}\min}_{\beta\in\R^d}\biggl(R^{(q)}_n(\beta
)+\lambda
(n,d,x)\frac{\|\beta\|_{\ell_1}^q}{n\epsilon^2}\biggr).
\]
Then, with probability greater than $1-12\exp(-x)$, the $L_q$-risk of
$\hat\beta_n$ satisfies
\[
R^{(q)}(\hat\beta_n)\leq\inf_{\beta\in\R^d}\biggl((1+2\epsilon
)R^{(q)}(\beta)+\eta(n,d,x)\frac{(1+\|\beta\|_{\ell
_1}^q)}{n\epsilon
^2}\biggr),
\]
where $\eta(n,d,x)=c_1 K(d)^q(\log n)^{(4q-2)/q}(\log d)^2(x+\log n)$.
\end{Theo}

Procedures based on the $\ell_1$-norm as a regularizing or constraint
function have been studied extensively in the last few years. We only
mention a small fraction of this very extensive body of work \cite
{MR2533469,MR2312149,MR2382644,MR2555200,MR2500227,MR2386087,MR2278363,MR2488351,MR1379242,MR2396809,MR2543687,MR2279469}.
In fact, it is almost impossible to make a proper comparison even with
the results mentioned in this partial list. Some of these results are
close enough in nature to Theorem~\ref{TheoC} to allow a comparison. In
particular, in \cite{MR2240689}, the authors prove that with high
probability, the LASSO satisfies an exact oracle inequality with a
residual term $\sim\|\beta\|_{\ell_1}/\sqrt{n}$ up to logarithm
factors, under tail assumptions on $Y$ and $X$. In~\cite{MR2312149},
upper bounds on the risks $\E[\langle X,\hat\beta_n-\beta_0 \rangle
^2]$ and
$\|\hat\beta_n-\beta_0\|_{\ell_1}$ were obtained for a weighted
LASSO $\hat\beta_n$ when $\E(Y|X)=\langle X,\beta_0
\rangle$ for
$\beta
_0$ with short support. Exact oracle inequalities for RERM using an
entropy-based criterion or on an $\ell_p$ criterion (with $p$ close to
$1$) were obtained in \cite{MR2509076,MR2500227} for any convex and
regular loss function and with fast rates. Similar bounds were obtained
in \cite{MR2396809} for a~RERM using a weighted $\ell_1$-criterion. In
\cite{MR2533469} it is shown that the LASSO and Dantzig estimators
\cite{MR2382644} satisfy oracle inequalities in the deterministic design
setup and under the REC condition. In fact, in most of these results
the authors obtained exact oracle inequalities with an optimal residual
term of $|\mathrm{Supp}(\beta_0)|(\log d)/n$, which is clearly better than
the rate $\|\beta\|_{\ell_1}^2/n$ obtained in Theorem \ref{TheoC} for the
quadratic loss and in the same context.

However, it is important to note that all these exact oracle
inequalities were obtained under an assumption that is similar in
nature to the Restricted Isometry Property (RIP), whereas in Theorem \ref{TheoC}
one does not need that kind of assumption on the design. Although it
seems strange that it is possible to obtain fast rates without RIP
there is nothing magical here. In fact, the isomorphic argument used to
prove Theorem \ref{TheoB} (and thus Theorem \ref{TheoC}) shows that the\vadjust{\goodbreak} random operator
$\beta\in\R^d \rightarrow n^{-1/2}\sum_{i=1}^n (Y_i-\langle
X_i,\beta\rangle)
e_i\in\R^n$ satisfies some sort of an RIP, which actually coincides
with the RIP property in the noise-free case $Y=\langle X,\beta_0
\rangle$ for an
isotropic design. This indicates that RIP is not the key property in
establishing oracle inequalities for the prediction risk, but rather,
the ``isomorphic profile'' of the problem at hand, which takes into
account the structure of the class of functions.

Finally, a word about notation. Throughout, we denote absolute
constants or constants that depend on other parameters by $c$, $C$, $c_1$,
$c_2$, etc. (and, of course, we will specify when a constant is
absolute and when it depends on other parameters). The values of these
constants may change from line to line. The notation $x\sim y$ (resp.,
$x\lesssim y$) means that there exist absolute constants $0<c<C$ such
that $cy\leq x\leq Cy$ (resp., $x\leq Cy$). If $b>0$ is a~parameter,
then $x\lesssim_b y$ means that $x\leq C(b) y $ for some constant
$C(b)$ depending only on $b$. We denote by $\ell_p^d$ the space $\R^d$
endowed with the $\ell_p$ norm $\|x\|_{\ell_p^d}=(\sum
_{j}|x_j|^p)^{1/p}$. The unit ball there is denoted by $B_p^d$ and
the unit Euclidean sphere in $\R^d$ is $S^{d-1}$.


\section{Preliminaries to the proofs}
\label{secpreliminaries}
In this section we obtain a general bound on $\E\|P-P_n\|_{(\ell
_F)_\lambda}$ for the $L_q$-loss when $q\geq2$, and show that a
Bernstein-type condition is satisfied under weak assumption on the loss
function.
\subsection{Isomorphic properties of the loss class}
\label{secisomorphy-class-function} The \textit{isomorphic property} of
a functions class measures the ``level'' at which empirical means and
actual means are equivalent. The notion was introduced in this context
in~\cite{MR2240689}. Although it is not a necessary feature of this
method, if one wishes the isomorphic property to hold with exponential
probability, one can use a~high probability deviation bound on the
supremum of the localized process. A~standard way (though not the only
way, or even the optimal way!) of obtaining such a~result is through of
Talagrand concentration inequality \cite{MR1258865} applied to
localizations of the function class, combined with a good control of
the variance in terms of the expectation (a Bernstein-type condition).
When applied to an excess loss class, this argument leads to exact
oracle inequalities; see, for example,
\cite{Mendelson08regularizationin,BMN}. Here we are interested in
nonexact oracle inequality, and thus, we will study the isomorphic
properties of the loss class. To make the presentation simpler, we are
not dealing with a fully ``unbounded theory'' like in
\cite{Shahar-Gaussian}, but rather that the class has an envelope
function which is bounded in~$\psi _1$, and we follow the path of
\cite{Mendelson08regularizationin}, in which one obtains the desired
high probability bounds using Talagrand's concentration theorem. Since
we would like to avoid the assumption that the class consists of
uniformly bounded functions, an important part of our analysis is the
following~$\psi_1$ version of Talagrand's inequality \cite{MR2424985}.
%
\begin{Theorem}\label{theoadamcjak}
There exists an absolute constant $K>0$ for which the following holds.
Let $Z_1,\ldots,Z_n$ be $n$ i.i.d. random variables with values in a~space $\cZ$,
and let $G$ be a countable class of real-valued measurable\vadjust{\goodbreak}
functions defined on $\cZ$. For every $x>0$ and $\alpha>0$, with
probability greater than $1-4\exp(-x)$,
\[
\|P-P_n\|_G\leq(1+\alpha)\E\|P-P_n\|_G+K\sigma(G)\sqrt
{\frac
{x}{n}}+K(1+\alpha^{-1})b_n(G)\frac{x}{n}.
\]
\end{Theorem}

Using the same truncation argument as in \cite{MR2424985}, it follows
that for every single function $g\in L_2(P)$ and every $\alpha,x>0$,
with probability greater than $1-4\exp(-x)$,
\[
P_n g\leq(1+\alpha)Pg+K\sqrt{\frac{x Pg^2}{n}}+K(1+\alpha
^{-1})\frac
{b_n(g)x}{n}
\]
and, in particular, if there exists some $ B_n\geq0$ for which
$Pg^2\leq B_n Pg+B_n^2/n$, then for every $0<\alpha<1$ and $x>0$, with
probability greater than $1-4\exp(-x)$,
%
\begin{equation}\label{eqadamcjak-single-function-under-Bernstein}
P_n g\leq(1+2\alpha)Pg+K^\prime(1+\alpha^{-1})\bigl(b_n(g)+B_n\bigr)\frac{x+1}{n}.
\end{equation}

Theorem \ref{theoadamcjak} can be extended to classes $G$ satisfying
some separability property like condition (M) in \cite{MR2291502}. We
apply Theorem \ref{theoadamcjak} in this context and it will be
implicitly assumed that every time we use Theorem \ref{theoadamcjak},
this separability condition holds. In particular, Theorem \ref
{theoadamcjak} will be applied to the localized sets $V(\ell
_F)_\lambda
$ to get nonexact oracle inequalities for the ERM algorithm and to the
family $(V(\ell_{F_r})_\lambda)_{r\geq0}$ to get nonexact regularized
oracle inequalities for the RERM procedure.

Observe that Theorem \ref{theoadamcjak} requires that the envelope
function $\sup_{g\in G}|g|$ is sub-exponential, but since $\|{\max_{1
\leq i \leq n} X_i}\|_{\psi_1} \lesssim\|X\|_{\psi_1}\log n$ it follows
that $b_n(\ell_F)$ is not much larger than $\|{\sup_{g \in G} g(X)}\|
_{\psi_1}$. However,\vspace*{1pt} this condition can be a major drawback. For
instance, if the set $G$ consists of linear functions indexed by the
Euclidean sphere $\cS^{d-1}$, and $X$ is the standard Gaussian measure
on $\R^d$, the resulting envelope function is bounded in $\psi_1(\mu)$,
but its norm is of the order of $\sqrt{d}$. In Theorem \ref{TheoC}, we bypass
this\vspace*{2pt} obstacle by assuming that $\|Y\|_{\psi_q},\|\|X\|_{\ell _\infty
^d}\|_{\psi_q}\leq K(d)$. This assumption is far better suited
for situations in which the indexing class is small---like localized
subsets of~$B_1^d$ that appear naturally in LASSO type results.
%
\begin{Theorem}
\label{corisomorphy}
Let $F$ be a functions class and assume that there exists $B_n\geq0$
such that for every $f\in F$, $ P \ell_f^2\leq B_n P\ell_f+B_n^2/n$. If
$0<\epsilon<1/2$ and $\lambda_\epsilon^*>0$ satisfy that
\[
\E\|P_n - P\|_{V(\ell_F)_{\lambda^*_\epsilon}} \leq(\epsilon
/4)\lambda
^*_\epsilon,
\]
then for every $x>0$, with probability larger than $1 -4 e^{-x}$, for
every $f \in F$
\[
P \ell_f \leq(1+2\epsilon)
P_n \ell_{f} + \rho_n(x),
\]
where, for $K$ the constant appearing in Theorem \ref{theoadamcjak},
\[
\rho_n(x)= \max\biggl(\lambda^*_\epsilon,\frac{(4Kb_n(\ell_F)+(6K)^2
B_n/\epsilon)(x+1)}{n\epsilon}\biggr).
\]
\end{Theorem}
\begin{pf}
The proof follows the ideas from \cite{MR2240689}. Fix $\lambda>0$ and
$x>0$, and note that by Theorem \ref{theoadamcjak}, with probability
larger than
$1-4\exp(-x)$,
%
\begin{eqnarray}\label{eqAdamcjak}
\|P - P_n\|_{V(\ell_F)_\lambda} &\leq&2
\E\|P-P_n\|_{V(\ell_F)_{\lambda}} + K \sigma(V(\ell
_F)_{\lambda})
\sqrt{\frac x n} \nonumber\\[-8pt]\\[-8pt]
&&{}+ K b_n(V(\ell_F)_{\lambda})\frac x n.\nonumber
\end{eqnarray}
Clearly, we have $b_n(V(\ell_F)_\lambda)\leq b_n(\ell_F)$ and
\[
\sigma^2(V(\ell_F)_{\lambda}) = \sup\bigl(P(\alpha\ell_{f})^2\dvtx 0
\leq
\alpha\leq1, f \in F, P( \alpha\ell_{f}) \leq\lambda\bigr) \leq
B_n\lambda+B_n^2/n.
\]

Moreover, since $V(\ell_F)$ is star-shaped, $\lambda\geq0\rightarrow
\phi(\lambda)= \E\|P-P_n\|_{V(\ell_F)_\lambda}/\lambda$ is
nonincreasing, and since $\phi(\lambda_\epsilon^*)\leq\epsilon/8$ and
$\rho_n(x)\geq\lambda^*_\epsilon$, then
\[
\E\|P-P_n\|_{V(\ell_F)_{\rho_n(x)}}\leq(\epsilon/4)\rho_n(x).
\]
Combined with (\ref{eqAdamcjak}), there
exists an event $\Omega_0(x)$ of probability greater than
$1-4\exp(-x)$, and on $\Omega_0(x)$,
\begin{eqnarray*}
\|P-P_n\|_{V(\ell_F)_{\rho_n(x)}} & \leq &
(\epsilon/2)\rho_n(x) + K
\sqrt{\frac{(B_n\rho_n(x)+B_n^2/n) x}{n}} + K\frac{b_n(\ell_F)x}{n}
\\
& \leq & \epsilon\rho_n(x).
\end{eqnarray*}
Hence, on $\Omega_0(x)$, if $g\in V(\ell_F)$ satisfies that $P g\leq
\rho
_n(x)$, then $|Pg - P_ng|\leq\epsilon\rho_n(x)$. Moreover, if
$P\ell_{f}=\beta>\rho_n(x)$, then
$g=\rho_n(x)\ell_{f}/\beta\in V(\ell_F)_{\rho_n(x)}$; hence
$|Pg-P_ng|\leq
\epsilon\rho_n(x)$, and so $(1-\epsilon)P\ell_{f}\leq P_n\ell
_{f}\leq
(1+\epsilon)P\ell_{f}$.
\end{pf}

\subsection{The Bernstein condition of loss functions classes}
\label{secBernstein-condition-loss-function-class}
In Theorem \ref{TheoA}, the desired concentration properties (and thus the fast
rates in Theorem \ref{TheoC}) rely on a Bernstein-type condition, that for every
$f \in F$,
%
\begin{equation}
\label{eqBernstein-assumption}
P\ell_f^2\leq B_n P\ell_f+B_n^2/n.
\end{equation}
%
%
%

Assumption (\ref{eqBernstein-assumption}) is trivially satisfied when
the loss functions are positive and uniformly bounded: if $0\leq\ell
_f\leq B$, then $P\ell_f^2\leq BP\ell_f$. It also turns out that
(\ref
{eqBernstein-assumption}) does not require any ``global'' structural
assumption on $F$ and is trivially verified if class members have
sub-exponential tails.
%
\begin{Lemma}
\label{lemBernstein-condition-under-psi1}
Let $X$ be a nonnegative subexponential random variable. Then for
every $z\geq1$,
\[
\E X^2\leq\log(ez)\|X\|_{\psi_1}\E X+\frac{(4+6\log
^2(ez)\|X\|_{\psi_1}^2)}{ez}.
\]
\end{Lemma}
\begin{pf} Fix $\theta> 0$, and note that
%
\begin{eqnarray}\label{eqlemma-bern-main}\quad
\E X^2\one_{X\geq\theta}&=&\int_0^\infty2t\Pro[X\one_{X\geq
\theta
}\geq
t]\,dt=\theta^2\Pro[X\geq\theta]+ 2\int_\theta^\infty
t\Pro
[X\geq t]\,dt\nonumber\\
&\leq& 2\theta^2\exp(-\theta/\|X\|_{\psi_1})+4\int
_\theta
^\infty t\exp(-t/\|X\|_{\psi_1})\,dt\\
&\leq&
(2\theta^2+4\theta\|X\|_{\psi_1}+4)\exp(-\theta
/\|X\|_{\psi_1}).\nonumber
\end{eqnarray}
Since $X\geq0$, it follows from (\ref{eqlemma-bern-main}) that, for
any $\theta>0$,
\begin{eqnarray*}
\E X^2&\leq&\E X^2\one_{X\leq\theta}+\E X^2\one_{X\geq\theta}\\
&\leq&\theta\E X+(2\theta^2+4\theta\|X\|_{\psi_1}+4
)\exp
(-\theta/\|X\|_{\psi_1}).
\end{eqnarray*}
The result follows for $\theta=\|X\|_{\psi_1}\log(ez)$.
\end{pf}

In particular, if $\ell_f \geq0$ and $\|\ell_f\|_{\psi_1}\leq
D$ for
some $D\geq1$, then for every $n\geq1$,
\[
\E\ell_f^2\leq(c_0D\log(en))\E\ell_f+\frac{
(c_0D\log
(en))^2}{n}.
\]


\subsection{\texorpdfstring{Upper bounds on $\E\|P-P_n\|_{V(\ell_F)_\lambda}$}{Upper bounds on E||P-P_n||_{V(l_F)_lambda}}}
\label{secUP-expectation}

Let $H$ be the loss class associated with $F$ for the ERM or with a
class $F_r$ for some $r\geq0$ for the RERM. The next step is to obtain
bounds on the fixed point of the localized process, that is, for some
$c_0<1$, to find a small $\lambda^*$ for which
\[
\E\|P-P_n\|_{V(H)_{\lambda^*}}\leq c_0\lambda^*.
\]

Note that the complexity of the star-shaped hull $V(H)$ is not far from
the one of $H$ itself. Actually, a bound on the expectation of the
supremum of the empirical process indexed by $V(H)_\lambda$ will follow
from one on $H_\mu$ for different levels $\mu\in\{2^i\lambda\dvtx i\in
\N\}$.
This follows from the peeling argument of \cite{BMN}:
that $ V(H)_\lambda\subset\bigcup_{i=0}^\infty\{\theta h\dvtx 0\leq
\theta
\leq2^{-i}, h\in H, \E h \leq2^{i+1}\lambda\}$.
Therefore, setting $H_\mu=\{h\in H\dvtx \E h\leq\mu\}$, for all $\mu>0$ and
$R^*=\inf_{h \in H} \E h$,
%
\begin{equation}\label{eqpeeling-shahar}
\E\|P-P_n\|_{V(H)_{\lambda}}
\leq\sum_{\{i\dvtx 2^{i+1}\lambda\geq R^*\}} 2^{-i}\E\|P-P_n\|
_{H_{2^{i+1}\lambda}},
\end{equation}
because if $2^{i+1}\lambda< R^*$, then the sets $H_{2^{i+1}\lambda}$
are empty. Thus, it remains to bound $\E\|P-P_n\|_{H_{\mu}}$ for any
$\mu>0$.

Let us mention that a naive attempt to control these empirical
processes using a contraction argument is likely to fail, and will
result in slow rates even in very simple cases (e.g., a regression
model with a bounded design). We refer to
\cite{MR1857312,shahar-psi1,MR2373017} for more details.


The bounds obtained below on $\E\|P-P_n\|_{H_\mu}$ are expressed in
terms of a~random metric complexity of $H$, which is based on the
structure of a~typical coordinate\vadjust{\goodbreak} projection $P_\sigma H$. These random
sets are defined for every sample $\sigma=(X_1,\ldots,X_n)$ by
\[
P_\sigma H=\{(f(X_1),\ldots,f(X_n))\dvtx f\in H\}.
\]

The complexity of these random sets will be measured via a metric
invariant, called the $\gamma_2$-functional, introduced by Talagrand as
a part of the generic chaining mechanism.
%
\begin{Definition}[(\cite{Talagrand05})]\label{defgamma2}
Let $(T,d)$ be a semi-metric space. An admissible sequence of $T$ is a
sequence\vspace*{1pt} $(T_s)_{s\in\N}$ of subsets of $T$ such that $|T_0|\leq1$ and
$|T_s|\leq2^{2^s}$ for any $s\geq1$. We define
\[
\gamma_2(T,d)=\inf_{(T_s)_{s\in\N}}\sup_{t\in T}\sum
_{s=0}^{\infty
}2^{s/2}d(t,T_s),
\]
where the infimum is taken over all admissible sequences $(T_s)_{s\in
\N
}$ of $T$.
\end{Definition}

We refer the reader to \cite{Talagrand05} for an extensive survey on
chaining methods and on the $\gamma_2$-functionals. In particular, one
can bound the $\gamma_2$-functional using an entropy integral
%
\begin{equation}
\label{eqDudley-gamma2}
\gamma_2(T,d)\lesssim\int_{0}^{\operatorname{diam}(T,d)} \sqrt{\log
N(T,d,\epsilon)}\,d\epsilon,
\end{equation}
where $ N(T,d,\epsilon)$ is the minimal number of balls of radius
$\epsilon$ with respect to the metric $d$ needed to cover $T$, and
$\operatorname{diam}(T,d)$ is the diameter of the metric space $(T,d)$.

We will use the $\gamma_2$-functional to state our theoretical bounds
because there are examples in which $\gamma_2(T,d)$ is significantly
smaller than the corresponding entropy integral. However, in all our
concrete applications we will use the bound (\ref{eqDudley-gamma2})
since the computation of those is much simpler, the gap is at most
logarithmic and the purpose of this note is not to obtain the optimal
estimates but to show that the residual terms in exact and nonexact
oracle inequalities could be very different.

Now, we turn to some concrete examples where $H$ is the loss functions
class in the regression model with respect to the $L_q$-loss.

Let\vspace*{-1pt} $q\geq2$ and set the $L_q$-loss function of $f$ to be $\ell
^{(q)}_f(x,y)= |y-f(x)|^q$. In this case, the $L_q$-loss functions
class localized at some level $\mu$ is $(\ell_F^{(q)})_\mu=\{\ell
_f^{(q)}\dvtx f\in F, \E\ell_f^{(q)}\leq\mu\}$.\vspace*{1pt}

The following result is a combination of a truncation argument and
Rudelson's $L_\infty^n$ method. To formulate it, set
$M=\|{\sup_{\ell\in(\ell_F^{(q)})_\mu}}|\ell|\|_{\psi_1}$, for
any\vspace*{-1pt} \mbox{$A
\subset\R^d$}, let $\widetilde{A}=A \cup-A$, and if $F^{(\mu)}=\{
f\in
F\dvtx P\ell_f^{(q)}\leq\mu\}$, put $U_n=\E\gamma_2^2(\widetilde
{P_\sigma
F}^{(\mu)} ,\ell_\infty^n)$.

\begin{Proposition}\label{propMark-learning}
For every $q\geq2$, there exists a constant $c_0$ depending only on
$q$ for which the following holds. If $F$ is a class of functions,\vadjust{\goodbreak} then
for any $\mu>0$:
\begin{longlist}[(2)]
\item[(1)] if $q=2$, then $\E\|P-P_n\|_{(\ell^{(q)}_F)_\mu}\leq c_0
\max[\sqrt{\mu\frac{U_n}{n}},\frac{U_n}{n}]$,
\item[(2)] if $q>2$, then $\E\|P-P_n\|_{(\ell^{(q)}_F)_\mu}$ is
upper bounded by
\[
c_0\max\Biggl[\sqrt{\mu\frac{U_n}{n}} \sqrt{(M\log n
)^{(q-2)/q}},\frac{U_n}{n} (M\log n)^{(q-2)/q},\frac{M\log
n}{n}\Biggr].
\]
\end{longlist}
%
\end{Proposition}
\begin{pf} Let $\phi(h)=\operatorname{sign}(h)\min(|h|,\theta)$ where $\theta
>0$ is a threshold to be fixed later. For $f\in F$, set
$h_f(x,y)=y-f(x)$, let \mbox{$H_\mu=\{h_f\dvtx f\in F, \E|h_f|^q\leq\mu\}$}, and
note that $|h|^q=|\phi(h)|^q+(|h|^q-\theta^q)\one_{|h|\geq
\theta
}$. Thus,
\begin{eqnarray*}
&&
\E\|P-P_n\|_{(\ell_F^{(q)})_\mu}\\
&&\qquad=\E\sup
_{h\in H_\mu
}\bigl|(P_n-P)(|h|^q)\bigr|\\
&&\qquad\leq\E\sup_{h\in H_\mu}\bigl|(P_n-P)(|\phi(h)|^q)\bigr|+\E
\sup
_{h\in H_\mu} P_n|h|^q\one_{|h|\geq\theta}\\
&&\qquad\quad{}+\sup_{h\in H_\mu
}P|h|^q\one
_{|h|\geq\theta}\\
&&\qquad\leq\E\sup_{h\in H_\mu}\bigl|(P_n-P)(|\phi(h)|^q)\bigr|+2\E
\Bigl(\sup_{h\in H_\mu} |h|^q\one_{|h|\geq\theta}\Bigr).
\end{eqnarray*}

To upper bound the truncated part of the process, consider the
empirical diameter $D_n=\sup_{h\in H_\mu}(P_n|\phi
(h)|^{2q-2}
)^{{1}/({2q-2})}$.
By the Zin{\'e}--Ginn symmetrization theorem \cite{vanderVaartWellner}
and the upper bound on a Rademacher process by a Gaussian one,
\[
\E\sup_{h\in H_\mu}\bigl|(P_n-P)(|\phi(h)|^q)\bigr|\leq\frac{c_0}{\sqrt
{n}}\E\E
_g\sup_{h\in H_\mu}\Biggl|\frac{1}{\sqrt{n}}\sum_{i=1}^ng_i|\phi
(h)(X_i,Y_i)|^q\Biggr|,
\]
where $g_1,\ldots,g_n$ are $n$ independent standard random variables
and $\E_g$ denotes the expectation with respect to those variables. For
a fixed sample $(X_i,Y_i)_{i=1}^n$, let $(Z(h))_{h\in H_\mu}$ be the
Gaussian process defined by $Z(h)= n^{-1/2}\sum_{i=1}^ng_i|\phi
(h)(X_i$, $Y_i)|^q$. If $f,g\in F$, then
\begin{eqnarray*}
&&\E_g\bigl(Z(h_f)-Z(h_g)\bigr)^2\\
&&\qquad=\frac{1}{n}\sum_{i=1}^n\bigl(|\phi
(h_f)(X_i,Y_i)|^q-|\phi(h_g)(X_i,Y_i)|^q\bigr)^2\\
&&\qquad\leq\frac{1}{n}\sum_{i=1}^nq^2|f(X_i)-g(X_i)|^2\max(|\phi
(h_f)(X_i,Y_i)|,|\phi(h_g)(X_i,Y_i)|)^{2q-2}\\
&&\qquad\leq 2q^2\max_{1\leq i\leq n}\bigl(f(X_i)-g(X_i)\bigr)^2 D_n^{2q-2},
\end{eqnarray*}
where we have used that $||\phi(u)|^q-|\phi(v)|^q|\leq
q|u-v|\max(|\phi(u)|,|\phi(v)|)^{q-1}$ for every $u,v\in
\R$.
By a standard chaining argument it follows that
%
\begin{equation}\label{eqgeneric-chaining-truncation-2}\qquad
\E_g\sup_{f\in F^{(\mu)}}\Biggl|\frac{1}{\sqrt{n}}\sum
_{i=1}^ng_i|\phi
(h_f)(X_i,Y_i)|^q\Biggr|\leq c_1 q\gamma_2\bigl(\widetilde{P_\sigma F}{}^{(\mu
)},\ell_\infty^n\bigr) D_n^{q-1}
\end{equation}
and thus, $\E\sup_{h\in H_\mu}|(P_n-P)(|\phi(h)|^q)|\leq c_2 q\sqrt{
\frac{\E\gamma_2^2(\widetilde{P_\sigma F}{}^{(\mu)},\ell_\infty
^n)}{n}}\sqrt{\E D_n^{2q-2}}$.

A bound on the diameter follows from (\ref
{eqgeneric-chaining-truncation-2}) and the contraction principle,
\begin{eqnarray*}
\E D_n^{2q-2} &\leq& \E\sup_{h\in H_\mu}\bigl|(P_n-P)(|\phi
(h)|^{2q-2})\bigr|+\sup
_{h\in H_\mu} P|\phi(h)|^{2q-2}\\
&\leq& \frac{c_2 q\theta^{q-2}}{\sqrt{n}}\E_g\sup_{h\in H_\mu
}\Biggl|\frac
{1}{\sqrt{n}}\sum_{i=1}^ng_i|\phi(h)(X_i,Y_i)|^q\Biggr|+\theta
^{q-2}\mu\\
&\leq& c_2q\theta^{q-2}\sqrt{\frac{U_n \E D_n^{2q-2}}{n}}+\theta
^{q-2}\mu,
\end{eqnarray*}
implying that $ \E D_n^{2q-2}\leq c_3\max(q^2\theta
^{2q-4}U_n/n,\theta^{q-2}\mu)$
and so
%
\begin{equation}
\label{eqbound-truncation-part}\qquad
\E\sup_{h\in H_\mu}\bigl|(P_n-P)(|\phi(h)|^q)\bigr|\leq c_4 q\max\Biggl(\frac
{qU_n\theta^{q-2}}{n},\sqrt{\frac{U_n\theta^{q-2}\mu}{n}}\Biggr).
\end{equation}
Next, observe that for $q=2$, the right-hand side in (\ref
{eqbound-truncation-part}) does not depend on the truncation level
$\theta$, and thus one may take $\theta$ arbitrarily large, leading to
the desired result.

For $q\neq2$, consider the unbounded part of the process. Since the
envelope function of $H_\mu$ exhibits a subexponential decay, then
\begin{eqnarray*}
\E\Bigl(\sup_{h\in H_\mu} |h|^q\one_{|h|\geq\theta}\Bigr)&=&\int
_0^\infty
\Pro\Bigl[\sup_{h\in H_\mu} |h|^q\one_{|h|\geq\theta}\geq t
\Bigr]\,dt\\
&=&\theta^q\Pro\Bigl[\sup_{h\in H_\mu}|h|\geq\theta\Bigr]+\int
_{\theta
^q}^\infty\Pro\Bigl[\sup_{h\in H_\mu}|h|^q\geq t\Bigr]\,dt\\
&\leq& 2\theta^q\exp(-\theta^q/M)+2M\exp(-\theta
^q/M).
\end{eqnarray*}
The result follows by taking $\theta^q =M \log n$.
\end{pf}



\section{\texorpdfstring{Proof of Theorem \protect\ref{TheoA}}{Proof of Theorem A}}
\label{secproof-TheoremA}
In this section, we will present the proof of Theorem~\ref{TheoA}, which follows
the same ideas as \cite{BMN,MR2240689} for the excess loss.
%
\begin{Lemma}\label{theosmall-lemma-ERM} There exists an absolute
constant $c_0>0$ for which the following holds. Let $F$ be a class of
functions, and assume that there is some $B_n$ such that for every
$f\in F$, $P\ell_f^2\leq B_n P\ell_f+B_n^2/n$. For $x>0$ and
$0<\epsilon
<1/2$, consider an event $\Omega_0(x)$ on which for every $f\in F$,
\[
R(f)\leq(1+2\epsilon)R_n(f)+\rho_n(x),\vadjust{\goodbreak}
\]
where $\rho_n(\cdot)$ is some fixed increasing function. Then, with
probability greater than $\Pro(\Omega_0(x))-4\exp(-x)$,
\[
R(\hat f^{\mathrm{ERM}}_n)\leq(1+3\epsilon)\inf_{f\in F}\biggl( R(f)+c_0\frac
{(b_n(\ell_f)+B_n)(x+1)}{n\epsilon}\biggr)+\rho_n(x).
\]
\end{Lemma}
\begin{pf} Fix $x>0$, let $K^\prime$ be the constant introduced in
(\ref{eqadamcjak-single-function-under-Bernstein}), consider
\[
f^* \in\mathop{\operatorname{Arg}\min}_{f\in F}\biggl(R(f)+15K^\prime\frac{ (b_n(\ell
_f)+B_n)(x+1)}{n\epsilon}\biggr)
\]
and without loss of generality one assume that the infimum is achieved.
By~(\ref{eqadamcjak-single-function-under-Bernstein}) [for $\alpha
=(\epsilon/2)/(1+2\epsilon)$], the event $\Omega^*(x)$ on which
\[
R_n(f^*)\leq\frac{1+3\epsilon}{1+2\epsilon} R(f^*)+5K^\prime\frac
{(b_n(\ell_{f^*})+B_n)(x+1)}{n\epsilon}
\]
has probability greater than $1-4\exp(-x)$.
Hence,
\[
-(1+3\epsilon)R(f^*)\leq-(1+2\epsilon)R_n(f^*)+15K^\prime\frac
{(b_n(\ell_{f^*})+B_n)(x+1)}{n\epsilon}
\]
and on $\Omega_0(x)\cap\Omega^*(x)$, every $f$ in $F$ satisfies that
\begin{eqnarray*}
R(f)-(1+3\epsilon) R(f^*)&\leq&(1+2\epsilon)\bigl(R_n(f)-R_n(f^*)
\bigr)+\rho_n(x)\\
&&{}+15K^\prime\frac{(b_n(\ell_{f^*})+B_n)(x+1)}{n\epsilon}.
\end{eqnarray*}
Since $R_n(\ERM)-R_n(f^*)\leq0$, then
\[
R(\ERM)\leq(1+3\epsilon) R(f^*)+15K^\prime\frac{(b_n(\ell
_{f^*})+B_n)(x+1)}{n\epsilon}+\rho_n(x),
\]
and the claim now follows from the choice of $f^*$.
\end{pf}
\begin{pf*}{Proof of Theorem \ref{TheoA}}
Let $x>0$, $0<\epsilon<1/2$, and put
\[
\rho_n(x)=\max\biggl(\lambda^*_\epsilon,\frac{(
(6K/\epsilon
)^2B_n+(4K/\epsilon)b_n(\ell_{F}))(x+1)}{n}\biggr).
\]
By Theorem \ref{corisomorphy}, the event $\Omega_0(x)$, on which every
$f \in F$ satisfies that
\[
R(f)\leq(1+2\epsilon)R_n(f)+\rho_n(x),
\]
has probability greater than $1-4\exp(-x)$.
Now, the result follows from Lem\-ma~\ref{theosmall-lemma-ERM}.

The remark following Theorem \ref{TheoA}, that if $\ell$ is
nonnegative, then $\ell_F$ satisfies a Bernstein-type condition with
$B_n \sim\operatorname{diam}(\ell_F,\psi_1)\log(en)$ follows from
Lem\-ma~\ref {lemBernstein-condition-under-psi1}.
\end{pf*}

\section{\texorpdfstring{Proof of Theorem \protect\ref{TheoB}}{Proof of Theorem B}}
\label{secoracle-ineq-RERM}

Although the proof of Theorem \ref{TheoB} seems rather technical, the idea behind
it is rather simple. First, one needs to find a~``trivial'' bound on
$\operatorname{crit}(\hat f^{\mathrm{RERM}}_n)$, giving preliminary information on where
one must look for the RERM function (this is the role played by the\vadjust{\goodbreak}
function~$\alpha_n$). Then, one combines peeling and fixed point
arguments to identify the exact location of the RERM.

Note that for $F=\bigcup_{r\geq0}F_r$, we have $\operatorname{crit}(f)=\infty$ for
all $f\in\cF\setminus F$. Therefore, without loss of generality, we
can replace the set $\cF$ by $F$ in both the definition of the RERM in
(\ref{eqRERM}) and in the nonexact regularized oracle inequality of
Theorem~\ref{TheoB}.

We begin with the following rough estimate on the criterion of the
RERM. In the case where there is a\vspace*{1pt} trivial bound $\operatorname{crit}(f)\leq
C_n$, for all $f\in F$ then it follows that for any $0<\epsilon<1/2$
and $x>0$, $\operatorname{crit}(\hat f^{\mathrm{RERM}}_n)\leq C_n=\alpha_n(\epsilon,x)$.
Turning to the second case stated in Assumption \ref{assthm-B}, recall
that $r\rightarrow\lambda_\epsilon^*(r)$ tends to infinity with $r$ and
there exists $K_1>0$ such that for every $(r,x)\in\R_+\times\R_+^*$,
$2\rho_n(r,x)\leq\rho_n(K_1(r+1),x)$. Hence, for every $x>0$ and
$0<\epsilon< 1/2$, we set $\alpha_n$ to satisfy that
\begin{eqnarray*}
&&\alpha_n(\epsilon,x)\geq\max\bigl[K_1\bigl(\operatorname{crit}(f_0)+2\bigr),\\
&&\hspace*{72pt}(\lambda_\epsilon^*)^{-1}\bigl((1+2\epsilon)\bigl(3R(f_0)+2K^\prime
\bigl(b_n(\ell
_{f_0})+B_n(\operatorname{crit}(f_0))\bigr)\\
&&\hspace*{254pt}{}\times\bigl((x+1)/n\bigr)\bigr)\bigr)\bigr],
\end{eqnarray*}
where $f_0$ is any fixed function in $F$ (e.g., when $0\in F$,
one may take $f_0=0$), and $(\lambda_\epsilon^*)^{-1}$ is the
generalized inverse function of $\lambda_\epsilon^*$. In this case, we
prove the following high probability bound on $\operatorname{crit}(\hat f^{\mathrm{RERM}}_n)$.
%
\begin{Lemma}\label{lempen-of-RERM-small}Assume that $r\rightarrow
\lambda_\epsilon^*(r)$ tends to infinity when $r$ tends to infinity and
that there exists $K_1>0$ such that for every $(r,x)\in\R_+\times\R
_+^*$, $2\rho_n(r,x)\leq\rho_n(K_1(r+1),x)$. Then, under the
assumptions of Theorem \ref{TheoB}, for every $x>0$ and $0<\epsilon< 1/2$, with
probability greater than $1-4\exp(-x)$, $ \operatorname{crit}(\hat
f^{\mathrm{RERM}}_n)\leq\alpha_n(\epsilon,x)$.
\end{Lemma}
\begin{pf}
By the definition of $\hat f^{\mathrm{RERM}}_n$,
\begin{eqnarray*}
&&R_n(\RERM)+\frac{2}{1+2\epsilon}\rho_n\bigl(\operatorname{crit}(\RERM
)+1,x+\log
\alpha_n(\epsilon,x)\bigr)\\
&&\qquad\leq R_n(f_0)+\frac{2}{1+2\epsilon}\rho_n\bigl(\operatorname{crit}(f_0)+1,x+\log
\alpha_n(\epsilon,x)\bigr).
\end{eqnarray*}
Since $\ell$ is nonnegative, then $R_n(\RERM)\geq0$, and thus
\begin{eqnarray*}
&&\rho_n\bigl(\operatorname{crit}(\RERM)+1,x+\log\alpha_n(\epsilon,x)\bigr)\\
&&\qquad\leq(1+2\epsilon)R_n(f_0)/2+\rho_n\bigl(\operatorname{crit}(f_0)+1,x+\log
\alpha
_n(\epsilon,x)\bigr)\\
&&\qquad\leq\max\bigl((1+2\epsilon)R_n(f_0),2\rho_n\bigl(\operatorname{crit}(f_0)+1,x+\log\alpha_n(\epsilon,x)\bigr)\bigr).
\end{eqnarray*}
Since $\rho_n(r,x) \geq\lambda^*_\epsilon(r)$, for all $r\geq0$, one
of the following two situations occurs: either
\[
\lambda_\epsilon^*(\operatorname{crit}(\RERM))\leq(1+2\epsilon)R_n(f_0)
\]
or, noting that for every $(r,x)\in\R_+\times\R_+^*$, $2\rho
_n(r,x)\leq
\rho_n(K_1(r+1),x)$, then
\begin{eqnarray*}
&&\rho_n\bigl(\operatorname{crit}(\RERM)+1,x+\log\alpha_n(\epsilon,x)
\bigr)\\
&&\qquad\leq
2\rho_n\bigl(\operatorname{crit}(f_0)+1,x+\log\alpha_n(\epsilon,x)\bigr)\\
&&\qquad\leq\rho_n\bigl(K_1\bigl(\operatorname{crit}(f_0)+2\bigr),x+\log\alpha_n(\epsilon
,x)\bigr),
\end{eqnarray*}
and since $\rho_n$ is monotone in $r$ then
$\operatorname{crit}(\RERM)\leq K_1(\operatorname{crit}(f_0)+2)$.

Hence, in both cases
%
\begin{equation}\label{eqlemma1-eq1}\qquad
\operatorname{crit}(\RERM)\leq\max\bigl((\lambda_\epsilon^*)^{-1}
\bigl((1+2\epsilon)R_n(f_0)\bigr),K_1\bigl(\operatorname{crit}(f_0)+2\bigr)\bigr).
\end{equation}

On the other hand, according to (\ref
{eqadamcjak-single-function-under-Bernstein}), with probability
greater than $1-4\exp(-x)$, $R_n(f_0)\leq3R(f_0)+2K^\prime(b_n(\ell
_{f_0})+B_n(\operatorname{crit}(f_0)))(x+1)/n$. The result follows by plugging
the last inequality in (\ref{eqlemma1-eq1}) and since $\lambda
_\epsilon
$ is nondecreasing.
\end{pf}

The next step is to find an ``isomorphic'' result for $\RERM$. The idea
is to divide the set given by the trivial estimate on $\operatorname{crit}(\RERM
)$ into level sets and analyze each piece separately.
%
\begin{Lemma} \label{lemrisk-RERM}
Under the assumptions of Theorem \ref{TheoB},
for every $x>0$, with probability greater than $1-8\exp(-x)$,
\[
R(\RERM)\leq(1+2\epsilon)R_n(\RERM)+\rho_n\bigl(\operatorname{crit}(\RERM
)+1,x+\log\alpha_n(\epsilon,x)\bigr).
\]
\end{Lemma}
\begin{pf} Let $\Omega_0(x)$ be the event
\[
\frac{R(\RERM)-R_n(\RERM)}{2\epsilon R_n(\RERM)+\rho_n(\operatorname{crit}(\RERM)+1,x+\log\alpha_n(\epsilon,x))}\geq1,
\]
and we will show that this event has the desired small probability.

Clearly,
\[
\Pro[\Omega_0(x)]\leq\Pro[\Omega_0(x)\cap\{\operatorname{crit}(\RERM
)\leq
\alpha_n(\epsilon,x)\}]+\Pro[\operatorname{crit}(\RERM)>\alpha
_n(\epsilon
,x)],
\]
and by Lemma \ref{lempen-of-RERM-small}, $\Pro[\operatorname{crit}(\RERM
)>\alpha_n(\epsilon,x)]\leq4\exp(-x)$ in the second case of
Assumption \ref{assthm-B} or $\Pro[\operatorname{crit}(\RERM)>\alpha
_n(\epsilon,x)]=0$ when there is a trivial bound on the criterion.
Therefore, in any case, we have $ \Pro[\operatorname{crit}(\RERM)>\alpha
_n(\epsilon,x)]\leq4\exp(-x)$.

Recall that $F_i=\{f \in F\dvtx\operatorname{crit}(f) \leq i\}$, for all $i\in\N$,
and since $\rho_n$ is monotone in $r$, then
\begin{eqnarray*}
&&\Pro[\Omega_0(x)\cap\{\operatorname{crit}(\RERM)\leq\alpha
_n(\epsilon,x)\}
]\\
&&\qquad\leq\sum_{i=0}^{\lfloor\alpha_n(\epsilon,x)\rfloor}\Pro
[\Omega
_0(x)\cap\{i\leq\operatorname{crit}(\RERM)\leq i+1\}]\\
&&\qquad\leq\sum_{i=0}^{\lfloor\alpha_n(\epsilon,x)\rfloor}\Pro
\bigl[\exists
f\in F_{i+1}\dvtx\\
&&\hspace*{78.5pt} R(f)\geq(1+2\epsilon)R_n(f)+\rho_n\bigl(i+1,x+\log\alpha
_n(\epsilon,x)\bigr)\bigr].
\end{eqnarray*}

By Theorem \ref{corisomorphy}, for every $t>0$ and $i\in\N$, with
probability greater than $1-4\exp(-t)$, for every $f\in F_{i+1}$,
$P\ell
_f\leq(1+2\epsilon)P_n\ell_f+\rho_n(i+1,t)$. In particular,
\begin{eqnarray*}
&&\Pro\bigl[\exists f\in F_{i+1}\dvtx R(f)\geq(1+2\epsilon)R_n(f)+\rho
_n\bigl(i+1,x+\log\alpha_n(\epsilon,x)\bigr)\bigr]\\
&&\qquad\leq4\exp\bigl(-\bigl(x+\log\alpha_n(\epsilon,x)\bigr)\bigr).
\end{eqnarray*}
Hence, the claim follows, since
\begin{eqnarray*}
&&\Pro[\Omega_0(x)\cap\{\operatorname{crit}(\RERM)
\leq\alpha_n(\epsilon,x)\}]\\
&&\qquad\leq\sum_{i=0}^{\lfloor\alpha_n(\epsilon,x)\rfloor}4\exp
\bigl(-\bigl(x+\log
\alpha_n(\epsilon,x)\bigr)\bigr)\leq4\exp(-x).
\end{eqnarray*}
\upqed\end{pf}
\begin{pf*}{Proof of Theorem \ref{TheoB}}
Let $x>0$ and $0<\epsilon<1$. Without loss of generality, we assume
that, for the constant $K^\prime$ defined in (\ref
{eqadamcjak-single-function-under-Bernstein}), there exists $f^*\in F$
minimizing the function
\begin{eqnarray*}
f\in F
&\longrightarrow& (1+3\epsilon)R(f)+\rho_n\bigl(\operatorname{crit}(f)+1,x+\log
\alpha_n(\epsilon,x)\bigr)\\
&&{}+6K^\prime\frac{(b_n(\ell_f)+B_n(\operatorname{crit}(f))(x+1)}{\epsilon n}.
\end{eqnarray*}

Let $\Omega^*(x)$ be the event on which
\[
R_n(f^*)\leq\frac{1+3\epsilon}{1+2\epsilon}R(f^*)\\
+K^\prime\frac
{(b_n(\ell_{f^*})+B_n(\operatorname{crit}(f^*)))(x+1)}{n}\biggl(\frac
{1+3\epsilon
}{\epsilon}\biggr).
\]
Since $f^*\in F_{\operatorname{crit}(f^*)}$, then $P\ell_{f^*}^2\leq B_n(\operatorname{crit}(f^*))P\ell_{f^*}+B_n^2(\operatorname{crit}(f^*))/n$, and by~(\ref
{eqadamcjak-single-function-under-Bernstein}) [applied with $\alpha
=\epsilon/(1+2\epsilon)$], $\Pro(\Omega^*(x))\geq1-4\exp(-x)$.

Consider the event $\Omega_0(x)$, on which
\begin{eqnarray*}
R(\RERM)&\leq&(1+2\epsilon)R_n(\RERM)\\
&&{}+\rho_n\bigl(\operatorname{crit}(\RERM
)+1,x+\log\alpha_n(\epsilon,x)\bigr)
\end{eqnarray*}
and observe that by Lemma \ref{lemrisk-RERM}, $\Pro[\Omega
_0(x)]\geq
1-8\exp(-x)$. Therefore, on $\Omega_0(x)\cap\Omega^*(x)$, we have
\begin{eqnarray*}
&&R(\RERM)+\rho_n\bigl(\operatorname{crit}(\hat f^{\mathrm{RERM}}_n)+1,x+\log\alpha
_n(\epsilon,x)\bigr)\\
&&\qquad\quad{}-(1+3\epsilon)R(f^*)\\
&&\qquad\leq (1+2\epsilon)\bigl(R_n(\RERM)-R_n(f^*)\bigr)\\
&&\qquad\quad{}+2\rho_n\bigl(\operatorname{crit}(\hat f^{\mathrm{RERM}}_n)+1,x+\log\alpha_n(\epsilon
,x)\bigr)\\
&&\qquad\quad{}+6K^\prime\frac{(b_n(\ell_{f^*})+B_n(\operatorname{crit}(f^*)))(x+1)}{\epsilon n}\\
&&\qquad\leq (1+2\epsilon)\biggl(R_n(\RERM)+\frac{2}{1+2\epsilon}\rho
_n\bigl(\operatorname{crit}(\hat f_n^{\mathrm{RERM}})+1,x+\log\alpha_n(\epsilon,x)\bigr)\\
&&\hspace*{71pt}\qquad\quad{}-R_n(f^*)-\frac{2}{1+2\epsilon}\rho_n\bigl(\operatorname{crit}(
f^*)+1,x+\log
\alpha_n(\epsilon,x)\bigr)\biggr)\\
&&\qquad\quad{}+2\rho_n\bigl(\operatorname{crit}(f^*)+1,x+\log\alpha_n(\epsilon,x)
\bigr)\\
&&\qquad\quad{}+6K^{\prime}\frac{(b_n(\ell_{f^*})+B_n(\operatorname{crit}(f^*)))(x+1)}{\epsilon n}\\
&&\qquad\leq 2\rho_n\bigl(\operatorname{crit}(f^*)+1,x+\log\alpha_n(\epsilon
,x)
\bigr)\\
&&\qquad\quad{}+6K^\prime\frac{(b_n(\ell_{f^*})+B_n(\operatorname{crit}(f^*)))(x+1)}{\epsilon
n},
\end{eqnarray*}
where the last inequality follows from the definition of $\hat
f^{\mathrm{RERM}}_n$. Hence, by the choice of $f^*$, it follows that on $\Omega
_1(x)\cap\Omega^*(x)$,
\begin{eqnarray*}
&&R(\RERM)+\rho_n\bigl(\operatorname{crit}(\hat f^{\mathrm{RERM}}_n)+1,x+\log\alpha
_n(\epsilon,x)\bigr)\\
&&\qquad\leq(1+3\epsilon)R(f^*)+2\rho_n\bigl(\operatorname{crit}(f^*)+1,x+\log
\alpha
_n(\epsilon,x)\bigr)\\
&&\qquad\quad{} +6K^\prime\frac{(b_n(\ell_{f^*})+B(\operatorname{crit}(f^*)))(x+1)}{\epsilon n}\\
&&\qquad= \inf_{f\in F}\biggl((1+3\epsilon)R(f)+ 2 \rho_n\bigl(\operatorname{crit}(f)+1,x+\log\alpha_n(\epsilon,x)\bigr)\\
&&\hspace*{85pt}\qquad\quad{} +6K^\prime\frac{(b_n(\ell_{f})+B(\operatorname{crit}(f)))(x+1)}{\epsilon
n}\biggr).
\end{eqnarray*}
\upqed\end{pf*}

\section{\texorpdfstring{Proofs of Theorem \protect\ref{TheoC}}{Proofs of Theorem C}}
\label{secproof-theo-C}
Theorem \ref{TheoC} follows from a direct application of Theorem \ref{TheoB}, by estimating
the specific function $\rho_n$ and the ``Bernstein function'' $B_n(r)$.

Consider the family of models $(F_r)_{r\geq0}$ associated with the
$\ell
_1$-criterion $F_r=\{f_\beta\dvtx\|\beta\|_1\leq r\}$,
where $f_\beta(x)=\langle x,\beta\rangle$ is a linear functional on
$\R^d$.
%
\begin{Lemma}\label{lemcomplexity-Elastic-Net}
There exists an absolute constant $c_0$ for which the following holds.
For every $\mu$ and $r\geq0$, and every $\sigma=(X_1,\ldots,X_n)$,
\[
\gamma_2(\widetilde{P_\sigma F}_r,\ell_\infty^n)\leq c_0r
\Bigl({\max
_{1\leq i\leq n}}\|X_i\|_{\ell_\infty^d}\Bigr)(\log d)\log
\biggl(\frac
{\sqrt{n}}{\log d}\biggr).
\]
Moreover, if $\|\|X\|_{\ell_\infty^d}\|_{\psi_2}\leq K(d)$, then
\[
(\E\gamma_2^2(\widetilde{P_\sigma F}_r,\ell_\infty^n)
)^{1/2}\leq c_0 rK(d)(\log n)^{3/2} (\log d).
\]
\end{Lemma}

The proof of the first part of the claim is rather standard and has
appeared in one form or another in several places; for example, see
\cite{BMN}. It follows from~(\ref{eqDudley-gamma2}) and Maurey's
empirical method (cf. \cite{MR810669,MR659306}). The second part is an
immediate corollary of the first one.

\begin{pf*}{Proof of Theorem \ref{TheoC}}
Observe that for every $\beta\in rB_1^d$,
\begin{eqnarray*}
\bigl\||Y-\langle X,\beta\rangle|^q\bigr\|_{\psi_1}&=&\|Y-\langle X,\beta
\rangle\|_{\psi
_q}^q\leq(\|Y\|_{\psi_q}+\|\langle X,\beta\rangle \|_{\psi
_q})^q\\
&\leq&(\|Y\|_{\psi_q}+\|\beta\|_1\|\|X\|_\infty \|_{\psi
_q})^q\leq(K(d))^q(1+r)^q.
\end{eqnarray*}
Hence, by Lemma \ref{lemBernstein-condition-under-psi1}, one may take
$B_n(r)=c_0(2K(d))^q(1+r)^q\log(en)$.

Next, the $\psi_1$-norm of the envelope of the class $F_r$ satisfies
$\|{\sup_{\beta\in rB_1^d}}|Y-\langle X,\beta\rangle|^q\|
_{\psi
_1}\leq
(K(d))^q(1+r)^q$,
and by (\ref{eqpeeling-shahar}), Proposition \ref{propMark-learning}
and Lemma \ref{lemcomplexity-Elastic-Net}, for every $\lambda>0$,
\begin{eqnarray*}
&&\E\|P-P_n\|_{V(\ell^{(q)}_{F_r})_\lambda}\\
&&\qquad\leq\sum
_{i=0}^\infty
2^{-i}\E\|P-P_n\|_{(\ell^{(q)}_{F_r})_{2^{i+1}\lambda}}\\
&&\qquad\leq c_0\sum_{i=0}^\infty2^{-i}\max\Biggl(\sqrt{2^{i+1}\lambda
}\sqrt
{\frac{r^2(1+r)^{q-2}h(n,d)}{n}},\\
&&\qquad\quad\hspace*{71pt}\frac{r^2(1+r)^{q-2}h(n,d)}{n},\frac{K(d)^q(1+r)^q (\log
n)}{n}\Biggr)\\
&&\qquad\leq c_1 \max\Biggl(\sqrt{\lambda}\sqrt{\frac
{(1+r)^qh(n,d)}{n}},\frac
{(1+r)^qh(n,d)}{n}\Biggr),
\end{eqnarray*}
where $h(n,d)=K(d)^q(\log n)^{(4q-2)/q}(\log d)^2$. Set $\lambda
_\epsilon^*(r)=c_2(1+r)^q h(n,d)/\allowbreak(n\epsilon^2)$
and observe that $\E\|P-P_n\|_{V(\ell^{(q)}_{F_r})_{\lambda
_\epsilon
^*(r)}}\leq(\epsilon/4)\lambda_\epsilon^*(r)$. Since
\[
b_n\bigl(\ell^{(q)}_{F_r}\bigr)=\Bigl\|\max_{1 \leq i \leq n} \sup_{f \in F_r} \ell
^{(q)}_f(X_i,Y_i)\Bigr\|_{\psi_1} \leq c_3(\log en)\Bigl\|{\sup_{f \in F_r} \ell
^{(q)}_f(X,Y)} \Bigr\|_{\psi_1},
\]
then one can take $\phi_n(r)=c_3K(d)^q(\log n)(1+r)^q$.
Thus
\[
\rho_n(r,x)=c_4\frac{h(n,d)(1+r^q)}{n\epsilon^2}(1+x)
\]
is a valid isomorphic function for this problem.
It is also easy to check that for \mbox{$f_0\equiv0$}, $\log\alpha
_n(\epsilon
,x)\leq c_5\log( \max(x,n)\|Y\|_{\psi_q}^q)$.
The result now follows by combining these estimates with Theorem \ref{TheoB}.
\end{pf*}

\section{Remarks on the differences between exact and nonexact oracle
inequalities}
\label{secComments}
The goal of this section is to describe the difference between the
analysis used in \cite{MR2240689} to obtain exact oracle inequalities
for the ERM,\vadjust{\goodbreak} and the one used in this note to establish nonexact
oracle inequalities for the ERM (Theorem \ref{TheoA}). Our aim is to indicate why
one may get faster rates for nonexact inequalities than for exact ones
for the same problem.

One should stress that this is not, by any means, a proof that it is
impossible to get exact oracle inequalities with fast rates (there are
in fact examples in which the ERM satisfies exact oracle inequalities
with fast rates: the linear aggregation problem, \cite{MR2329442}). It
is not even a proof that the localization method presented here is
sharp. A detailed study of the isomorphic method and oracle
inequalities for a general sub-Gaussian case (i.e., a~sub-exponential
squared loss), in the sense that the class $F$ has a bounded diameter
in $L_{\psi_2}$ rather than an envelope function, will be presented in
\cite{Shahar-Gaussian}.

However, we believe that this explanation will help to shed some light
on the differences between the two types of inequalities, and we refer
the reader to \cite{Shahar-Gaussian} for a more detailed and accurate analysis.

Our starting point is the following exact oracle inequality for ERM,
which is a mild modification of a result from \cite{MR2240689}. The
only difference is that it uses Adamczak's $\psi_1$ version of
Talagrand's concentration inequality for empirical processes, instead
of Massart's version.
%
\begin{Theorem}\label{theoExact-oracle-ineq-ERM}
There exists an absolute constant $c_0>0$ for which the
following holds. Let $F$ be a class of functions and assume that there
exists $B>0$ such that for every $ f\in F$, $P\cL_f^2\leq BP\cL_f$. Let
$\mu^*>0$ be such that $\E\|P_n - P\|_{V(\cL_F)_{\mu^*}} \leq\mu^*/8$,
and consider an increasing function $\rho_n$ which satisfies that, for
every $x>0$, $\rho_n(x)\geq\max(\mu^*,c_0(b_n(\cL
_F)+B)x/n)$.
Then, for every $x>0$, with probability greater than $1-8\exp(-x)$, the
risk of the ERM satisfies $R(\hat f_n^{\mathrm{ERM}})\leq\inf_{f\in
F}R(f)+\rho_n(x)$.
\end{Theorem}

Roughly put, and as indicated by the theorem, localization arguments
are based on two main components:

\begin{longlist}[(2)]
\item[(1)] A Bernstein-type condition, the essence of which is that it
allows one to ``translate'' localization with respect to the loss or the
excess loss to a~localization with respect to a~natural metric. In
particular this leads to the necessary control on the $\ell_2^n$
diameter of a random coordinate projection of the localized class.

\item[(2)] The fixed point of the empirical process indexed by the
localized star-shaped hull of the loss functions class (for nonexact
inequalities) or of the excess loss functions class (for exact ones).
\end{longlist}

Although the two components seem similar for the exact and nonexact
cases, they are very different. Indeed, for a nonexact oracle
inequality, the Bernstein type condition is almost trivially satisfied
and requires no special properties on the model/output couple
$(F,Y)$---as long as the functions involved have well behaved tails. As such,
it is an individual property of every class member; see Lem\-ma~\ref
{lemBernstein-condition-under-psi1}.\vadjust{\goodbreak}

On the other hand, the Bernstein condition required for the exact
oracle inequality is deeply connected to the geometry of the problem;
see, for example, \cite{MR2426759}. More accurately, when the target
$Y$ is far from the set of multiple minimizers of the risk, $N(F,\ell
,X)=\{Y\dvtx|\{f\in F\dvtx R(f)=\inf_{f\in F}R(f)\}|\geq2\}$, one can show that
a Bernstein condition holds for a large variety of loss function $\ell
$. However, when the target $Y$ gets closer to the set $N(F,\ell,X)$,
the Bernstein constant $B$ degenerates, and leads to rates slower than
$1/\sqrt{n}$ even if $F$ is a two functions class. Hence, the geometry
of the problem (the relative position of $Y$ and $F$) is very important
when trying to establish exact oracle inequalities, and the Bernstein
condition is truly a ``global'' property of $F$.

In particular, this explains the gap that we observed in the example
preceding the formulation of Theorem \ref{TheoA}. In that case, the class is a
finite set of functions and the set $N(F,\ell,X)$ is nonempty. Thus,
one can find a set~$F$ and a target $Y$ in a ``bad'' position, leading
to an excess loss class $\cL_F$ with a trivial Bernstein constant
(i.e., greater than $\sqrt{n}$). On the other hand, regardless of the
choice of $Y$, the Bernstein constant of $\ell_F$ is well behaved.

Let us mention that when the gap between exact and nonexact oracle
inequalities is only due to the Bernstein condition, it is likely that
both ERM and RERM will be suboptimal procedures \cite
{lbw96,MR2451042,LM2}. In particular, when slow rates are due to a lack
of convexity of $F$ (which is closely related to a bad Bernstein
constant of ${\cal L}_F$), one can consider procedures which ``improve
the geometry'' of the model (e.g., the ``starification'' method of
\cite{Audibert1} or the ``pre-selection-convexification'' method
in~\cite{LM1}).

The second aspect of the problem is the fixed point of the localized
empirical process. Although the complexity of the sets $\cL_F$ and
$\ell
_F$ seems similar from a metric point of view ($\cL_F$ is just a shift
of $\ell_F$) the localized star-shaped hull $(\cL_F)_\lambda$ and
$(\ell
_F)_\lambda$ are rather different. Since there are many ways of
bounding the empirical process indexed by these localized sets, let us
show the difference for one of the methods---based on the random
geometry of the classes, and for the sake of simplicity, we will only
consider the square loss. Using this method of analysis at hand, the
dominant term of the bound on $\E\|P-P_n\|_{V(\ell_F^{(2)})_\mu}$
(for the loss class) which was obtained in Proposition~\ref
{propMark-learning} is
%
\begin{equation}\label{eqmain-term-loss-functions-class}
\sqrt{\mu}\sqrt{\frac{\E\gamma_2^2(\widetilde{P_\sigma F}{}^{(\mu
)},\ell
_\infty^n)}{n}}.
\end{equation}
A similar bound was obtained for $\E\|P-P_n\|_{V(\cL_F)_\mu}$ in
\cite{Mendelson08regularizationin} and \cite{BMN}, in which the
dominant term is
%
\begin{equation}\label{eqmain-term-excess-loss-function-class}
\sqrt{\Bigl(\inf_{f\in F}R(f)+\mu\Bigr)}\sqrt{\frac{\E\gamma
_2^2(\widetilde{P_\sigma F}{}^{(\mu)},\ell_\infty^n)}{n}}.
\end{equation}
If this bound is sharp (and it is in many cases), and since
$R^*=\inf_{f\in F}R(f)$ is in general a nonzero constant, the fixed
point $\mu^*$ of Theorem \ref{theoExact-oracle-ineq-ERM} is of the
order of $\sqrt{\E\gamma_2^2(\widetilde{P_\sigma F}{}^{(\mu^*)},\ell
_\infty^n)/n}$ and thus leads to a rate decaying more slowly than
$1/\sqrt{n}$. In contrast, in the nonexact case one has $\lambda
_\epsilon^*\sim\E\gamma_2^2(\widetilde{P_\sigma F}{}^{(\lambda
_\epsilon
^*)},\ell_\infty^n)/n$ which is of the order of $1/n$ (up to
logarithmic factors) when the complexity $\E\gamma_2^2(\widetilde
{P_\sigma F},\ell_\infty^n)$ is ``reasonable.''

The reason for this gap comes from the observation that functions in
the star hull of $\ell_F$ whose expectation is smaller than $R^*$ are
only ``scaled down'' versions of functions from $\ell_F$. In fact, the
``complexity'' of the localized sets below the level of $R^*$ can
already be seen at the level $R^*$. Hence, the empirical process those
sets index (when scaled properly), becomes smaller with $\lambda$.

In contrast, because there are functions ${\cal L}_f$ that can have an
arbitrarily small expectation, the complexity of the localized subsets
of the star hull of~${\cal L}_F$ (normalized properly, of course) can
even increase as $\lambda$ decreases. This happens in very simple
situations; for example, even in regression relative to~$B_1^d$, if
$R^* \not= 0$, the complexity of the localized sets remains almost
stable and starts to decrease only at a very ``low'' level $\lambda$.
This is the reason for the phase transition in the error rate
(\mbox{$\sim$}$\max\{\sqrt{(\log d)/n}, d/n\}$) that one encounters in that problem.
The first term is due to the fact that the complexity of the localized
sets does not change as $\lambda$ decreases---up to some critical
level, while the second captures what happens when the localized sets
begin to ``shrink.''
A concrete example of this phenomenon is treated in the Supplementary
material \cite{supplementary} in the Convex aggregation context.


\begin{supplement}
\stitle{Applications to matrix completion, convex aggregation and model
selection}
\slink[doi]{10.1214/11-AOS965SUPP} 
\sdatatype{.pdf}
\sfilename{aos965\_supp.pdf}
\sdescription{In the supplementary file, we apply our main results to
the problem of matrix completion, convex aggregation and model
selection. The aim is to expose the fundamental differences between
exact and nonexact oracle inequalities on classical problems.}
\end{supplement}


\printaddresses


\begin{thebibliography}{47}

\bibitem{MR2424985}
\begin{barticle}[mr]
\bauthor{\bsnm{Adamczak},~\bfnm{Rados{\l}aw}\binits{R.}}
(\byear{2008}).
\btitle{A tail inequality for suprema of unbounded empirical processes with
  applications to {M}arkov chains}.
\bjournal{Electron. J. Probab.}
\bvolume{13}
\bpages{1000--1034}.
\bid{doi={10.1214/EJP.v13-521}, issn={1083-6489}, mr={2424985}}
\bptok{imsref}%
\end{barticle}
\endbibitem

\bibitem{Audibert1}
\begin{bmisc}[auto:STB|2012/03/21|07:41:58]
\bauthor{\bsnm{Audibert},~\bfnm{Jean-Yves}\binits{J.-Y.}}
(\byear{2007}).
\bhowpublished{No fast exponential deviation inequalities for the progressive
  mixture rule. Technical report, CERTIS.}
\bptok{imsref}%
\end{bmisc}
\endbibitem

\bibitem{MR1679028}
\begin{barticle}[mr]
\bauthor{\bsnm{Barron},~\bfnm{Andrew}\binits{A.}},
  \bauthor{\bsnm{Birg{\'e}},~\bfnm{Lucien}\binits{L.}} \AND
  \bauthor{\bsnm{Massart},~\bfnm{Pascal}\binits{P.}}
(\byear{1999}).
\btitle{Risk bounds for model selection via penalization}.
\bjournal{Probab. Theory Related Fields}
\bvolume{113}
\bpages{301--413}.
\bid{doi={10.1007/s004400050210}, issn={0178-8051}, mr={1679028}}
\bptok{imsref}%
\end{barticle}
\endbibitem

\bibitem{MR2240689}
\begin{barticle}[mr]
\bauthor{\bsnm{Bartlett},~\bfnm{Peter~L.}\binits{P.~L.}} \AND
  \bauthor{\bsnm{Mendelson},~\bfnm{Shahar}\binits{S.}}
(\byear{2006}).
\btitle{Empirical minimization}.
\bjournal{Probab. Theory Related Fields}
\bvolume{135}
\bpages{311--334}.
\bid{doi={10.1007/s00440-005-0462-3}, issn={0178-8051}, mr={2240689}}
\bptok{imsref}%
\end{barticle}
\endbibitem

\bibitem{BMN}
\begin{bmisc}[mr]
\bauthor{\bsnm{Bartlett},~\bfnm{Peter~L.}\binits{P.~L.}},
  \bauthor{\bsnm{Mendelson},~\bfnm{Shahar}\binits{S.}} \AND
    \bauthor{\bsnm{Neeman},~\bfnm{Joseph}\binits{J.}}
(\byear{2012}).
\bhowpublished{$\ell_1$-regularized linear regression: Persistence and oracle
inequalities.
\textit{Probab. Theory Related Fields}. To appear.}
\bptok{imsref}%
\end{bmisc}
\endbibitem

\bibitem{MR2533469}
\begin{barticle}[mr]
\bauthor{\bsnm{Bickel},~\bfnm{Peter~J.}\binits{P.~J.}},
  \bauthor{\bsnm{Ritov},~\bfnm{Ya'acov}\binits{Y.}} \AND
  \bauthor{\bsnm{Tsybakov},~\bfnm{Alexandre~B.}\binits{A.~B.}}
(\byear{2009}).
\btitle{Simultaneous analysis of lasso and {D}antzig selector}.
\bjournal{Ann. Statist.}
\bvolume{37}
\bpages{1705--1732}.
\bid{doi={10.1214/08-AOS620}, issn={0090-5364}, mr={2533469}}
\bptok{imsref}%
\end{barticle}
\endbibitem

\bibitem{MR2312149}
\begin{barticle}[mr]
\bauthor{\bsnm{Bunea},~\bfnm{Florentina}\binits{F.}},
  \bauthor{\bsnm{Tsybakov},~\bfnm{Alexandre}\binits{A.}} \AND
  \bauthor{\bsnm{Wegkamp},~\bfnm{Marten}\binits{M.}}
(\byear{2007}).
\btitle{Sparsity oracle inequalities for the {L}asso}.
\bjournal{Electron. J. Stat.}
\bvolume{1}
\bpages{169--194}.
\bid{doi={10.1214/07-EJS008}, issn={1935-7524}, mr={2312149}}
\bptok{imsref}%
\end{barticle}
\endbibitem

\bibitem{MR2382644}
\begin{barticle}[mr]
\bauthor{\bsnm{Candes},~\bfnm{Emmanuel}\binits{E.}} \AND
  \bauthor{\bsnm{Tao},~\bfnm{Terence}\binits{T.}}
(\byear{2007}).
\btitle{The {D}antzig selector: Statistical estimation when {$p$} is much
  larger than {$n$}}.
\bjournal{Ann. Statist.}
\bvolume{35}
\bpages{2313--2351}.
\bid{doi={10.1214/009053606000001523}, issn={0090-5364}, mr={2382644}}
\bptok{imsref}%
\end{barticle}
\endbibitem

\bibitem{MR810669}
\begin{barticle}[mr]
\bauthor{\bsnm{Carl},~\bfnm{Bernd}\binits{B.}}
(\byear{1985}).
\btitle{Inequalities of {B}ernstein--{J}ackson-type and the degree of
  compactness of operators in {B}anach spaces}.
\bjournal{Ann. Inst. Fourier (Grenoble)}
\bvolume{35}
\bpages{79--118}.
\bid{issn={0373-0956}, mr={0810669}}
\bptok{imsref}%
\end{barticle}
\endbibitem

\bibitem{MR2241189}
\begin{barticle}[mr]
\bauthor{\bsnm{Donoho},~\bfnm{David~L.}\binits{D.~L.}}
(\byear{2006}).
\btitle{Compressed sensing}.
\bjournal{IEEE Trans. Inform. Theory}
\bvolume{52}
\bpages{1289--1306}.
\bid{doi={10.1109/TIT.2006.871582}, issn={0018-9448}, mr={2241189}}
\bptok{imsref}%
\end{barticle}
\endbibitem

\bibitem{MR1857312}
\begin{bincollection}[mr]
\bauthor{\bsnm{Gin{\'e}},~\bfnm{Evarist}\binits{E.}},
  \bauthor{\bsnm{Lata{\l}a},~\bfnm{Rafa{\l}}\binits{R.}} \AND
  \bauthor{\bsnm{Zinn},~\bfnm{Joel}\binits{J.}}
(\byear{2000}).
\btitle{Exponential and moment inequalities for {$U$}-statistics}.
In \bbooktitle{High Dimensional Probability, {II} ({S}eattle, {WA}, 1999)}.
\bseries{Progress in Probability}
\bvolume{47}
\bpages{13--38}.
\bpublisher{Birkh\"auser}, \baddress{Boston, MA}.
\bid{mr={1857312}}
\bptok{imsref}%
\end{bincollection}
\endbibitem

\bibitem{MR2329442}
\begin{barticle}[mr]
\bauthor{\bsnm{Koltchinskii},~\bfnm{Vladimir}\binits{V.}}
(\byear{2006}).
\btitle{Local {R}ademacher complexities and oracle inequalities in risk
  minimization}.
\bjournal{Ann. Statist.}
\bvolume{34}
\bpages{2593--2656}.
\bid{doi={10.1214/009053606000001019}, issn={0090-5364}, mr={2329442}}
\bptok{imsref}%
\end{barticle}
\endbibitem

\bibitem{MR2555200}
\begin{barticle}[mr]
\bauthor{\bsnm{Koltchinskii},~\bfnm{Vladimir}\binits{V.}}
(\byear{2009}).
\btitle{The {D}antzig selector and sparsity oracle inequalities}.
\bjournal{Bernoulli}
\bvolume{15}
\bpages{799--828}.
\bid{doi={10.3150/09-BEJ187}, issn={1350-7265}, mr={2555200}}
\bptok{imsref}%
\end{barticle}
\endbibitem

\bibitem{MR2509076}
\begin{barticle}[mr]
\bauthor{\bsnm{Koltchinskii},~\bfnm{Vladimir}\binits{V.}}
(\byear{2009}).
\btitle{Sparse recovery in convex hulls via entropy penalization}.
\bjournal{Ann. Statist.}
\bvolume{37}
\bpages{1332--1359}.
\bid{doi={10.1214/08-AOS621}, issn={0090-5364}, mr={2509076}}
\bptok{imsref}%
\end{barticle}
\endbibitem

\bibitem{MR2500227}
\begin{barticle}[mr]
\bauthor{\bsnm{Koltchinskii},~\bfnm{Vladimir}\binits{V.}}
(\byear{2009}).
\btitle{Sparsity in penalized empirical risk minimization}.
\bjournal{Ann. Inst. Henri Poincar\'e Probab. Stat.}
\bvolume{45}
\bpages{7--57}.
\bid{doi={10.1214/07-AIHP146}, issn={0246-0203}, mr={2500227}}
\bptok{imsref}%
\end{barticle}
\endbibitem

\bibitem{MR1857339}
\begin{bincollection}[mr]
\bauthor{\bsnm{Koltchinskii},~\bfnm{Vladimir}\binits{V.}} \AND
  \bauthor{\bsnm{Panchenko},~\bfnm{Dmitriy}\binits{D.}}
(\byear{2000}).
\btitle{Rademacher processes and bounding the risk of function learning}.
In \bbooktitle{High Dimensional Probability, {II} ({S}eattle, {WA}, 1999)}.
\bseries{Progress in Probability}
\bvolume{47}
\bpages{443--457}.
\bpublisher{Birkh\"auser}, \baddress{Boston, MA}.
\bid{mr={1857339}}
\bptok{imsref}%
\end{bincollection}
\endbibitem

\bibitem{LM4}
\begin{bmisc}[auto:STB|2012/03/21|07:41:58]
\bauthor{\bsnm{Lecu{\'e}},~\bfnm{Guillaume}\binits{G.}} \AND
  \bauthor{\bsnm{Mendelson},~\bfnm{Shahar}\binits{S.}}
(\byear{2012}).
\bhowpublished{On the optimality of the aggregate with exponential weights for
  low temperature. \textit{Bernoulli}. To appear.}
\bptok{imsref}%
\end{bmisc}
\endbibitem

\bibitem{LM1}
\begin{barticle}[mr]
\bauthor{\bsnm{Lecu{\'e}},~\bfnm{Guillaume}\binits{G.}} \AND
  \bauthor{\bsnm{Mendelson},~\bfnm{Shahar}\binits{S.}}
(\byear{2009}).
\btitle{Aggregation via empirical risk minimization}.
\bjournal{Probab. Theory Related Fields}
\bvolume{145}
\bpages{591--613}.
\bid{doi={10.1007/s00440-008-0180-8}, issn={0178-8051}, mr={2529440}}
\bptnote{check year}%
\bptok{imsref}%
\end{barticle}
\endbibitem

\bibitem{LM2}
\begin{barticle}[mr]
\bauthor{\bsnm{Lecu{\'e}},~\bfnm{Guillaume}\binits{G.}} \AND
  \bauthor{\bsnm{Mendelson},~\bfnm{Shahar}\binits{S.}}
(\byear{2010}).
\btitle{Sharper lower bounds on the performance of the empirical risk
  minimization algorithm}.
\bjournal{Bernoulli}
\bvolume{16}
\bpages{605--613}.
\bid{doi={10.3150/09-BEJ225}, issn={1350-7265}, mr={2730641}}
\bptok{imsref}%
\end{barticle}
\endbibitem

\bibitem{supplementary}
\begin{bmisc}[auto:STB|2012/03/21|07:41:58]
\bauthor{\bsnm{Lecu{\'e}},~\bfnm{Guillaume}\binits{G.}} \AND
  \bauthor{\bsnm{Mendelson},~\bfnm{Shahar}\binits{S.}}
(\byear{2012}).
\bhowpublished{Supplement to ``General non-exact oracle inequalities for
  classes with a subexponential envelope.''
  DOI:\href{http://dx.doi.org/10.1214/11-AOS965SUPP}{10.1214/}
  \href{http://dx.doi.org/10.1214/11-AOS965SUPP}{11-AOS965SUPP}.}
\bptok{imsref}%
\end{bmisc}
\endbibitem

\bibitem{LT91}
\begin{bbook}[mr]
\bauthor{\bsnm{Ledoux},~\bfnm{Michel}\binits{M.}} \AND
  \bauthor{\bsnm{Talagrand},~\bfnm{Michel}\binits{M.}}
(\byear{1991}).
\btitle{Probability in {B}anach Spaces: Isoperimetry and Processes}.
\bseries{Ergebnisse der Mathematik und Ihrer Grenzgebiete (3) [Results in
  Mathematics and Related Areas (3)]}
\bvolume{23}.
\bpublisher{Springer}, \baddress{Berlin}.
\bid{mr={1102015}}
\bptok{imsref}%
\end{bbook}
\endbibitem

\bibitem{MR2386087}
\begin{barticle}[mr]
\bauthor{\bsnm{Lounici},~\bfnm{Karim}\binits{K.}}
(\byear{2008}).
\btitle{Sup-norm convergence rate and sign concentration property of {L}asso
  and {D}antzig estimators}.
\bjournal{Electron. J. Stat.}
\bvolume{2}
\bpages{90--102}.
\bid{doi={10.1214/08-EJS177}, issn={1935-7524}, mr={2386087}}
\bptok{imsref}%
\end{barticle}
\endbibitem

\bibitem{MR2319879}
\begin{bbook}[mr]
\bauthor{\bsnm{Massart},~\bfnm{Pascal}\binits{P.}}
(\byear{2007}).
\btitle{Concentration Inequalities and Model Selection}.
\bseries{Lecture Notes in Math.}
\bvolume{1896}.
\bpublisher{Springer}, \baddress{Berlin}.
\bid{mr={2319879}}
\bptok{imsref}%
\end{bbook}
\endbibitem

\bibitem{MR2291502}
\begin{barticle}[mr]
\bauthor{\bsnm{Massart},~\bfnm{Pascal}\binits{P.}} \AND
  \bauthor{\bsnm{N{\'e}d{\'e}lec},~\bfnm{{\'E}lodie}\binits{{\'E}.}}
(\byear{2006}).
\btitle{Risk bounds for statistical learning}.
\bjournal{Ann. Statist.}
\bvolume{34}
\bpages{2326--2366}.
\bid{doi={10.1214/009053606000000786}, issn={0090-5364}, mr={2291502}}
\bptok{imsref}%
\end{barticle}
\endbibitem

\bibitem{MR2278363}
\begin{barticle}[mr]
\bauthor{\bsnm{Meinshausen},~\bfnm{Nicolai}\binits{N.}} \AND
  \bauthor{\bsnm{B{\"u}hlmann},~\bfnm{Peter}\binits{P.}}
(\byear{2006}).
\btitle{High-dimensional graphs and variable selection with the lasso}.
\bjournal{Ann. Statist.}
\bvolume{34}
\bpages{1436--1462}.
\bid{doi={10.1214/009053606000000281}, issn={0090-5364}, mr={2278363}}
\bptok{imsref}%
\end{barticle}
\endbibitem

\bibitem{MR2488351}
\begin{barticle}[mr]
\bauthor{\bsnm{Meinshausen},~\bfnm{Nicolai}\binits{N.}} \AND
  \bauthor{\bsnm{Yu},~\bfnm{Bin}\binits{B.}}
(\byear{2009}).
\btitle{Lasso-type recovery of sparse representations for high-dimensional
  data}.
\bjournal{Ann. Statist.}
\bvolume{37}
\bpages{246--270}.
\bid{doi={10.1214/07-AOS582}, issn={0090-5364}, mr={2488351}}
\bptok{imsref}%
\end{barticle}
\endbibitem

\bibitem{Shahar-Gaussian}
\begin{bmisc}[auto:STB|2012/03/21|07:41:58]
\bauthor{\bsnm{Mendelson},~\bfnm{Shahar}\binits{S.}}
\bhowpublished{Oracle inequalities and the isomorphic method.
Technical report, Technion,
Israel Inst. Technology.}
\bptok{imsref}%
\end{bmisc}
\endbibitem

\bibitem{MR2451042}
\begin{barticle}[mr]
\bauthor{\bsnm{Mendelson},~\bfnm{Shahar}\binits{S.}}
(\byear{2008}).
\btitle{Lower bounds for the empirical minimization algorithm}.
\bjournal{IEEE Trans. Inform. Theory}
\bvolume{54}
\bpages{3797--3803}.
\bid{doi={10.1109/TIT.2008.926323}, issn={0018-9448}, mr={2451042}}
\bptok{imsref}%
\end{barticle}
\endbibitem

\bibitem{m08}
\begin{barticle}[mr]
\bauthor{\bsnm{Mendelson},~\bfnm{Shahar}\binits{S.}}
(\byear{2008}).
\btitle{Lower bounds for the empirical minimization algorithm}.
\bjournal{IEEE Trans. Inform. Theory}
\bvolume{54}
\bpages{3797--3803}.
\bid{doi={10.1109/TIT.2008.926323}, issn={0018-9448}, mr={2451042}}
\bptok{imsref}%
\end{barticle}
\endbibitem

\bibitem{MR2426759}
\begin{barticle}[mr]
\bauthor{\bsnm{Mendelson},~\bfnm{Shahar}\binits{S.}}
(\byear{2008}).
\btitle{Obtaining fast error rates in nonconvex situations}.
\bjournal{J. Complexity}
\bvolume{24}
\bpages{380--397}.
\bid{doi={10.1016/j.jco.2007.09.001}, issn={0885-064X}, mr={2426759}}
\bptok{imsref}%
\end{barticle}
\endbibitem

\bibitem{shahar-psi1}
\begin{barticle}[mr]
\bauthor{\bsnm{Mendelson},~\bfnm{Shahar}\binits{S.}}
(\byear{2010}).
\btitle{Empirical processes with a bounded {$\psi\sb 1$} diameter}.
\bjournal{Geom. Funct. Anal.}
\bvolume{20}
\bpages{988--1027}.
\bid{doi={10.1007/s00039-010-0084-5}, issn={1016-443X}, mr={2729283}}
\bptok{imsref}%
\end{barticle}
\endbibitem

\bibitem{Mendelson08regularizationin}
\begin{barticle}[mr]
\bauthor{\bsnm{Mendelson},~\bfnm{Shahar}\binits{S.}} \AND
  \bauthor{\bsnm{Neeman},~\bfnm{Joseph}\binits{J.}}
(\byear{2010}).
\btitle{Regularization in kernel learning}.
\bjournal{Ann. Statist.}
\bvolume{38}
\bpages{526--565}.
\bid{doi={10.1214/09-AOS728}, issn={0090-5364}, mr={2590050}}
\bptok{imsref}%
\end{barticle}
\endbibitem

\bibitem{MR2373017}
\begin{barticle}[mr]
\bauthor{\bsnm{Mendelson},~\bfnm{Shahar}\binits{S.}},
  \bauthor{\bsnm{Pajor},~\bfnm{Alain}\binits{A.}} \AND
  \bauthor{\bsnm{Tomczak-Jaegermann},~\bfnm{Nicole}\binits{N.}}
(\byear{2007}).
\btitle{Reconstruction and subgaussian operators in asymptotic geometric
  analysis}.
\bjournal{Geom. Funct. Anal.}
\bvolume{17}
\bpages{1248--1282}.
\bid{doi={10.1007/s00039-007-0618-7}, issn={1016-443X}, mr={2373017}}
\bptok{imsref}%
\end{barticle}
\endbibitem

\bibitem{Shahar-Paouris}
\begin{bmisc}[auto:STB|2012/03/21|07:41:58]
\bauthor{\bsnm{Mendelson},~\bfnm{Shahar}\binits{S.}} \AND
  \bauthor{\bsnm{Paouris},~\bfnm{Grigoris}\binits{G.}}
(\byear{2011}).
\bhowpublished{On the generic chaining and the smallest singular value of
  random matrices with heavy tails. Unpublished manuscript. Available at
  arXiv:\arxivurl{1108.3886}.}
\bptok{imsref}%
\end{bmisc}
\endbibitem

\bibitem{MR1320206}
\begin{barticle}[mr]
\bauthor{\bsnm{Natarajan},~\bfnm{B.~K.}\binits{B.~K.}}
(\byear{1995}).
\btitle{Sparse approximate solutions to linear systems}.
\bjournal{SIAM J. Comput.}
\bvolume{24}
\bpages{227--234}.
\bid{doi={10.1137/S0097539792240406}, issn={0097-5397}, mr={1320206}}
\bptok{imsref}%
\end{barticle}
\endbibitem

\bibitem{MR659306}
\begin{bincollection}[mr]
\bauthor{\bsnm{Pisier},~\bfnm{G.}\binits{G.}}
(\byear{1981}).
\btitle{Remarques sur un r\'esultat non publi\'e de {B}. {M}aurey}.
In \bbooktitle{Seminar on {F}unctional {A}nalysis, 1980--1981}
\bpages{Exp. No. V, 13}.
\bpublisher{\'Ecole Polytech.}, \baddress{Palaiseau}.
\bid{mr={0659306}}
\bptok{imsref}%
\end{bincollection}
\endbibitem

\bibitem{MR2450103}
\begin{bbook}[mr]
\bauthor{\bsnm{Steinwart},~\bfnm{Ingo}\binits{I.}} \AND
  \bauthor{\bsnm{Christmann},~\bfnm{Andreas}\binits{A.}}
(\byear{2008}).
\btitle{Support Vector Machines}.
\bpublisher{Springer}, \baddress{New York}.
\bid{mr={2450103}}
\bptok{imsref}%
\end{bbook}
\endbibitem

\bibitem{MR1258865}
\begin{barticle}[mr]
\bauthor{\bsnm{Talagrand},~\bfnm{M.}\binits{M.}}
(\byear{1994}).
\btitle{Sharper bounds for {G}aussian and empirical processes}.
\bjournal{Ann. Probab.}
\bvolume{22}
\bpages{28--76}.
\bid{issn={0091-1798}, mr={1258865}}
\bptok{imsref}%
\end{barticle}
\endbibitem

\bibitem{Talagrand05}
\begin{bbook}[mr]
\bauthor{\bsnm{Talagrand},~\bfnm{Michel}\binits{M.}}
(\byear{2005}).
\btitle{The Generic Chaining: Upper and Lower Bounds of Stochastic Processes}.
\bpublisher{Springer}, \baddress{Berlin}.
\bid{mr={2133757}}
\bptok{imsref}%
\end{bbook}
\endbibitem

\bibitem{MR1379242}
\begin{barticle}[mr]
\bauthor{\bsnm{Tibshirani},~\bfnm{Robert}\binits{R.}}
(\byear{1996}).
\btitle{Regression shrinkage and selection via the lasso}.
\bjournal{J. Roy. Statist. Soc. Ser. B}
\bvolume{58}
\bpages{267--288}.
\bid{issn={0035-9246}, mr={1379242}}
\bptok{imsref}%
\end{barticle}
\endbibitem

\bibitem{MR2396809}
\begin{barticle}[mr]
\bauthor{\bparticle{van~de} \bsnm{Geer},~\bfnm{Sara~A.}\binits{S.~A.}}
(\byear{2008}).
\btitle{High-dimensional generalized linear models and the lasso}.
\bjournal{Ann. Statist.}
\bvolume{36}
\bpages{614--645}.
\bid{doi={10.1214/009053607000000929}, issn={0090-5364}, mr={2396809}}
\bptok{imsref}%
\end{barticle}
\endbibitem

\bibitem{vanderVaartWellner}
\begin{bbook}[mr]
\bauthor{\bparticle{van~der} \bsnm{Vaart},~\bfnm{Aad~W.}\binits{A.~W.}} \AND
  \bauthor{\bsnm{Wellner},~\bfnm{Jon~A.}\binits{J.~A.}}
(\byear{1996}).
\btitle{Weak Convergence and Empirical Processes: With Applications to
  Statistics}.
\bpublisher{Springer}, \baddress{New York}.
\bid{mr={1385671}}
\bptok{imsref}%
\end{bbook}
\endbibitem

\bibitem{MR672244}
\begin{bbook}[mr]
\bauthor{\bsnm{Vapnik},~\bfnm{Vladimir}\binits{V.}}
(\byear{1982}).
\btitle{Estimation of Dependences Based on Empirical Data}.
\bpublisher{Springer}, \baddress{New York}.
\bid{mr={0672244}}
\bptok{imsref}%
\end{bbook}
\endbibitem

\bibitem{lbw96}
\begin{binproceedings}[auto:STB|2012/03/21|07:41:58]
\bauthor{\bsnm{Wee},~\bfnm{S.~Lee}\binits{S.~L.}},
  \bauthor{\bsnm{Bartlett},~\bfnm{Peter~L.}\binits{P.~L.}} \AND
  \bauthor{\bsnm{Williamson},~\bfnm{Robert~C.}\binits{R.~C.}}
(\byear{1996}).
\btitle{The importance of convexity in learning with squared loss}.
In \bbooktitle{Proceedings of the Ninth Annual Conference on Computational
  Learning Theory}
\bpages{140--146}.
\bpublisher{ACM Press}, \baddress{New York}.
\bptok{imsref}%
\end{binproceedings}
\endbibitem

\bibitem{MR2543687}
\begin{barticle}[mr]
\bauthor{\bsnm{Zhang},~\bfnm{Tong}\binits{T.}}
(\byear{2009}).
\btitle{Some sharp performance bounds for least squares regression with {$L\sb
  1$} regularization}.
\bjournal{Ann. Statist.}
\bvolume{37}
\bpages{2109--2144}.
\bid{doi={10.1214/08-AOS659}, issn={0090-5364}, mr={2543687}}
\bptok{imsref}%
\end{barticle}
\endbibitem

\bibitem{MR2279469}
\begin{barticle}[mr]
\bauthor{\bsnm{Zou},~\bfnm{Hui}\binits{H.}}
(\byear{2006}).
\btitle{The adaptive lasso and its oracle properties}.
\bjournal{J. Amer. Statist. Assoc.}
\bvolume{101}
\bpages{1418--1429}.
\bid{doi={10.1198/016214506000000735}, issn={0162-1459}, mr={2279469}}
\bptok{imsref}%
\end{barticle}
\endbibitem

\end{thebibliography}
\end{document}